\documentclass[12pt]{article}
\usepackage[margin=2.5cm]{geometry}
\usepackage{CJK}
\usepackage{latexsym,bm}
\usepackage{hyperref}
\usepackage{tabularx}
\usepackage{booktabs}
\usepackage{dcolumn}
\usepackage{amsfonts}
\usepackage{latexsym}
\usepackage{amsmath}
\usepackage{epsf}
\usepackage[section]{placeins}
\usepackage{mathabx}
\usepackage{graphicx}
\usepackage{subfigure}%(P229)
\usepackage{wrapfig}%(P227)
\usepackage{multicol}
\usepackage{verbatim}
\usepackage{mathrsfs}
\usepackage{multirow}
\usepackage{epstopdf}
\usepackage[subfigure]{graphfig}
\usepackage{stfloats}
\usepackage{float}
\usepackage{xcolor}
\usepackage{framed}
\usepackage{algorithm,algpseudocode,lipsum}
\usepackage{cite}
\floatname{algorithm}{Algorithm}

\usepackage[colorinlistoftodos]{todonotes} 

\numberwithin{equation}{section}

\newtheorem{assumption}{Assumption}
\newtheorem{lemma}{Lemma}
\newtheorem{theorem}{Theorem}[section]



\begin{document}

\title{Parallel in time partially explicit splitting scheme for high contrast linear multiscale diffusion problems}

 \author{Yating Wang \thanks{School of Mathematics and Statistic, Xi'an Jiaotong University, Xi'an, People's Republic of China.}
 \and Zhengya Yang \thanks{School of Mathematics and Statistic, Xi'an Jiaotong University, Xi'an, People's Republic of China.}
 \and Wing Tat Leung \thanks{Department of Mathematics, City University of Hong Kong, Hong Kong Special Administrative Region. }}

\maketitle

\begin{abstract}
Solving multiscale diffusion problems is often computationally expensive due to the spatial and temporal discretization challenges arising from high-contrast coefficients. To address this issue, a partially explicit temporal splitting scheme is proposed. By appropriately constructing multiscale spaces, the spatial multiscale property is effectively captured, and it has been demonstrated that the temporal step size is independent of the contrast. To enhance simulation speed, we propose a parallel algorithm for the multiscale flow problem that leverages the partially explicit temporal splitting scheme. The idea is first to evolve the partially explicit system using a coarse time step size, then correct the solution on each coarse time interval with a fine propagator, for which we consider the all-at-once solver. This procedure is then performed iteratively till convergence. We analyze the stability and convergence of the proposed algorithm. The numerical experiments demonstrate that the proposed algorithm achieves high numerical accuracy for high-contrast problems and converges in a relatively small number of iterations. The number of iterations stays stable as the number of coarse intervals increases, thus significantly improving computational efficiency through parallel processing.
\end{abstract}

\section{Introduction}\label{section:Introduction}
Numerous scientific problems and models exhibit multiscale properties, such as flow in heterogeneous porous media, the diffusion of pollutants in the atmosphere, turbulent transport in high Reynolds number flows, and so on. These models often involve significant variations in media properties, commonly referred to as high contrast. The presence of high contrast introduces stiffness to the system, which makes spatial and temporal discretization challenging for traditional numerical methods due to their high computational demands. Since it requires a small enough mesh size to capture multiscale features, and also a small enough time step to handle the stiffness arising from the high contrast.

\indent There have been many existing approaches in the literature to handle spatial multiscale problems, such as numerical homogenization (NH) \cite{reference1,reference2}, multiscale finite element methods (MsFEM) \cite{reference3,reference4,reference30}, generalized multiscale finite element methods (GMsFEM) \cite{reference5,reference6,reference7}, constraint energy minimizing GMsFEM (CEM-GMsFEM) \cite{reference8,reference9,reference25}, localized orthogonal decomposition (LOD) \cite{reference12,reference26} and nonlocal multi-continua method (NLMC) \cite{reference11,reference32,reference34}. Among which, CEM-GMsFEM is a multiscale finite element method used to effectively address multiscale problems with high-contrast parameters. It constructs multiscale basis functions by minimizing energy constraints, which can achieve contrast independent convergence rates. Based on CEM-GMsFEM, NLMC is proposed to construct local basis functions that automatically identify physical properties in each local region and provides non-local transmissibilities in the global formulation.

\indent For non-stationary multiscale problems, the high-contrast will lead to very small time steps when treating explicitly. The partially explicit temporal splitting scheme \cite{reference14} originates from the solution decomposition \cite{reference23} and splitting algorithms \cite{reference24}. The main idea of the method is to divide the solution space into two parts, the dominant part and the complementary part, such that the time step size is independent of the high-contrast. The method is successfully employed in solving wave equation, nonlinear diffusion, and time-fractional diffusion equations \cite{reference15,reference16,reference27,reference34,reference36} and extended to an adaptive algorithm\cite{reference34}. In this paper, we follow the concept developed in \cite{reference14} for linear equations.

\indent All the existing literature is based on the sequential method for solving the partially explicit temporal splitting scheme, which is easy to implement but might be inefficient when the temporal mesh partition is small enough and a long time simulation is needed. For this reason, we introduce a parareal algorithm to enhance computational efficiency. The parareal algorithm was proposed by Lions et al in \cite{reference13}. Its core idea is to divide the entire time interval into subintervals and compute simultaneously and independently on each subinterval. Numerous studies have investigated the analysis and applications of the parareal method, see \cite{reference17,reference19,reference20,reference21,reference31,reference43}. However, the existing literature is mainly based on Euler or Runge-Kutta method as the fine solver in each subinterval, which could be time-consuming for multiscale problems if one takes small time steps. To address this issue, we further introduce the waveform relaxation (WR) method \cite{reference38,reference39,reference40} via the diagonalization technique based on the all-at-once algorithm \cite{reference22, reference28,reference29} as the fine solver in the parareal framework \cite{reference35}. All-at-once algorithm is a global method that solves the problem over the entire time interval simultaneously instead of solving it step by step. It naturally fits for parallel computation and can significantly improve computational efficiency.

\indent The main contributions of this paper are as follows. \begin{itemize}
    \item The parareal algorithm for solving the partially explicit temporal splitting scheme is presented. The convergence of the proposed algorithm is shown. 
    \item The WR method via the diagonalization technique based on the all-at-once method is introduced into the parareal algorithm to speed up the computation for the fine propagator.
    \item The proposed algorithm can achieve high accuracy, and convergence can be reached with a small number of iterations. As the number of coarse interval (processors) increases, the number of iterations needed to achieve the error tolerance is quite stable, thus significantly saving the computational cost with parallel computation.
\end{itemize} 

\indent The rest of this paper is organized as follows. In Section \ref{section:Preliminaries}, we present preliminaries. In Section \ref{section:Multiscale space construction}, we give in detail the construction of multiscale spaces. The parareal all-at-once algorithm and the error estimate of full discretization are given in Section \ref{section:Parallel in time for partially explicit temporal splitting scheme}. Numerical experiments and conclusion are given in Section \ref{section:Numerical experiments} and 
Section \ref{section:Conclusion}, respectively.

\section{Preliminaries}\label{section:Preliminaries}
\indent In this paper, we consider the diffusion equation of the form
 \begin{align}\label{eqnAa}
 \left\{\begin{aligned}
 \frac{\partial u}{\partial t}-\nabla\cdot(\kappa\nabla u)&=f,\ \text{in}\ \Omega\times(0,T]\\
 u&=0,\ \text{on}\ \partial\Omega\times(0,T]\\
 u&=u_{0},\ \text{on}\ \partial\Omega\times\{0\}
 \end{aligned}\right.
 \end{align}
 where $\Omega$ is a bounded domain, $\kappa\in L^{\infty}(\Omega)$ is a high contrast parameter and $f\in L^{2}(0,T;L^{2}(\Omega))$ is the source term. 
  
 \indent We first present the fully-discretized problem for \eqref{eqnAa}, with the finite element method for the spatial discretization and backward Euler method for temporal discretization. Next, we derive the partially explicit temporal splitting scheme. Then we briefly introduce the framework of the parareal algorithm.
 
 \indent We now introduce some notations. Denote by $(u,v):=\displaystyle\int_{\Omega}uvd\Omega$ the inner product of $L^{2}(\Omega)$ whose norm is denoted by $\parallel\cdot\parallel_{\Omega}$. We use $H_{0}^{1}(\Omega)$ to denote the Sobolev spaces with zero boundary values. To simplify the notation, $\Omega$ may be dropped in the notations. We denote $C$ a generic positive constant independent of any function and of any discretization parameters.
 
 We write the problem \eqref{eqnAa} in the semi-discrete form: find $u(\cdot,t)\in H_{0}^{1}(\Omega)$ such that
 \begin{equation*}
 \left(\frac{\partial u}{\partial t},v\right)+a\left(u,v\right)=\left(f,v\right),\ \forall v\in H_{0}^{1}(\Omega),
 \end{equation*}
 where the bilinear form $a(\cdot,\cdot)$ is given by
 \begin{equation*}
 a(u,v)=\int_{\Omega}\kappa\nabla u\cdot\nabla v,
 \end{equation*}
 and define the energy norm $\parallel u\parallel_{a}=a(u,u)^{\frac{1}{2}}$.

 \indent Consider a coarse spatial partition $\mathcal{T}_{H}$ of $\Omega$ with mesh size $H$, we will construct multiscale basis functions on $\mathcal{T}_{H}$ and form a multiscale space $V_{H}$ which has good approximation power. For the approximation space $V_H\in H_{0}^{1}(\Omega)$, the semi-discretization in space leads to seeking $u_{H}(\cdot, t)\in V_{H}$ such that
 \begin{equation}\label{eqnAb}
 \begin{aligned}
 &\left(\frac{\partial u_{H}}{\partial t},v\right)+a\left(u_{H},v\right)=\left(f,v\right),\ \forall v\in V_H,\ 0<t\leq T,\\
 &u_{H}(0)=u_{H,0}.
 \end{aligned}
 \end{equation}
where $u_{H,0}$ is the projection of $u_{0}$ in $V_{H}$.

 Let $\Delta t$ be the time step size and $t_{n}=n\Delta t, n=0,1,\cdots,N, T=N\Delta t$. Then the full discretization with Backward Euler reads
 \begin{equation}\label{eqnAc}
 \left(\frac{u_{H}^{n+1}-u_{H}^{n}}{\Delta t},v\right)+a\left(u_{H}^{n+1},v\right)=\left(f^{n+1},v\right),\ \forall v\in V_{H},
 \end{equation}
 where $u_{H}^{n+1}\approx u(\cdot,t_{n+1})$. %It is well-known that the Backward Euler scheme is unconditionally stable. \\
 
 \subsection{Partially explicit temporal splitting scheme}
 \indent Now we introduce the partially explicit temporal splitting scheme \cite{reference14}. Assume that $V_{H}$ can be decomposed into two subspaces $V_{H,1}$ and $V_{H,2}$, that is 
 \begin{align*}
 V_{H} = V_{H,1}+V_{H,2}.
 \end{align*}
 Then a partially explicit temporal splitting scheme reads: finding $\{u_{n}\}_{n=1}^{N}\in V_{H,1}$ and $\{w_{n}\}_{n=1}^{N}\in V_{H,2}$ such that
 \begin{align}
 \label{eqnAd}\left(\frac{u_{n+1}-u_{n}}{\Delta t},v_{1}\right)+\left(\frac{w_{n}-w_{n-1}}{\Delta t},v_{1}\right)+a\left(u_{n+1}+w_{n},v_{1}\right)&=\left(f^{n},v_{1}\right),\\
 \label{eqnAe}\left(\frac{w_{n+1}-w_{n}}{\Delta t},v_{2}\right)+\left(\frac{u_{n}-u_{n-1}}{\Delta t},v_{2}\right)+a\left(u_{n+1}+w_{n},v_{2}\right)&=\left(f^{n},v_{2}\right),
 \end{align}
 $\forall v_{1}\in V_{H,1},\forall v_{2}\in V_{H,2}$. Initial conditions are projected onto corresponding subspaces. Thus, the solution at time step $n+1$ will be $u_{H}^{n+1}=u_{n+1}+w_{n+1}$.\\
 \indent It has been proved in \cite{reference14} that with a suitable choice of $V_{H,1}$ and $V_{H,2}$ the partially explicit temporal splitting scheme is stable. 
 
 \subsection{Parareal algorithm}
 \label{section:Parareal algorithm}
 \indent To enhance the computational efficiency of the partially explicit temporal splitting scheme \eqref{eqnAd}-\eqref{eqnAe}, we will introduce a temporal parallel algorithm. We first describe the basic flow of the parareal algorithm by considering the following initial value problem
 \begin{align}\label{eqnAf}
 \begin{aligned}
 &\frac{du}{dt}=F\left(t,u\right),\ t\in(0,T],\\
 &u(0)=u_{0}.
 \end{aligned}
 \end{align}
 We first divide $(0,T]$ into $N$ elements as described in the \eqref{eqnAc}. Let $\mathcal{F}$ be a fine solver that can achieve the desired accuracy but has a high computational cost. We also introduce a coarse solver $\mathcal{G}$, which have lower accuracy but provide results at a lower computational costs compared the fine solver. For example, one can use the same numerical scheme for both $\mathcal{F}$ and $\mathcal{G}$ but $\mathcal{F}$ using a small step size $\delta t$, while $\mathcal{G}$ is utilizing a bigger step size $\Delta t$ ($\Delta t\gg\delta t$).\\
\indent Denote by $\mathcal{F}_{\Delta t}(u,t_{n},t_{n+1})$ the result obtained by integrating $u$ from $t_{n}$ to $t_{n+1}$ using $\mathcal{F}_{\delta t}$ (fine solver with time step size $\delta t$) and $\mathcal{G}_{\Delta t}(u,t_{n},t_{n+1})$ denoting the similar integration forward in time using the coarse solver. At the zero iteration of the parareal method, we get $\{u_{n}^{0}\}_{n=0}^{N}$ using the coarse solver $\mathcal{G}$. Denote by $u_{n}^{k}$ the approximation for $u(t_{n})$ at the $k$th iteration. Then the solution of the $k+1$th iteration is obtained from the following formula: for all $0\leq n\leq N-1$ and $k=0,1,2,\cdots,$
\begin{equation}\label{eqnAg}
u_{n+1}^{k+1}=\mathcal{G}_{\Delta t}\left(u_{n}^{k+1},t_{n},t_{n+1}\right)+\mathcal{F}_{\Delta t}\left(u_{n}^{k},t_{n},t_{n+1}\right)-\mathcal{G}_{\Delta t}\left(u_{n}^{k},t_{n},t_{n+1}\right).
\end{equation}
 One can terminate the parareal algorithm if the maximum number of iterations is met ($k=n_{\text{max}}$), or if
\begin{align}\label{eqnAh}
\max_{1\leq n\leq N}\parallel u_{n}^{k}-u_{n}^{k-1}\parallel\leq\epsilon,
\end{align}
where $\epsilon$ is a given tolerance.

Notice that the evolution of the fine solver $\mathcal{F}$ only requires the initial value $u_{n}^{k}$, which depends on the previous iteration. Therefore, for each $k$, the $\mathcal{F}$ can be computed in parallel \cite{reference20,reference35}. \\
\indent The solution of \eqref{eqnAg} converges under suitable assumptions, i.e., $u_{n}^{k}\rightarrow u_{n}^{*}$, where $u_{n}^{*}$ is the solution obtained from $\mathcal{F}$ with the time step $\delta t$ throughout the temporal domain. 
\begin{figure}[H]
\centering
\begin{tikzpicture}[scale=0.45]
  \draw[step=0.5cm,gray,very thin] (-5,-5) grid (5,5);
  \draw[step=1cm,red,thin] (-5.5,-5.5) grid (5.5,5.5);
  \draw[step=2cm,blue,thick] (-5.5,-5.5) grid (5.5,5.5);
        \fill[olive!80!black, fill opacity=0.5] (-2,-2) rectangle (4,4); 
        \fill[yellow!90!orange, fill opacity=0.9] (0,0) rectangle (2,2); 
  % \draw[fill=yellow] (0,0) rectangle (2,2);
  \draw[ultra thick] (-2,-2) rectangle (4,4);

  \draw[fill=red] (3,3) rectangle (3.5,3.5);
  \node[font=\fontsize{8}{9.6}\selectfont\bfseries] at (3.98,4.4) {fine grid};
  \node[font=\fontsize{10}{12}\selectfont\bfseries] at (-3.15,-3.15) {$K_{i}^{+}$};
  \node[font=\fontsize{10}{12}\selectfont\bfseries] at (1,1) {$K_i$};

  \draw[->,thick] (3.9,4.2) -- (3.6,3.6);
  \draw[->,thick] (-3,-3) -- (-2.1,-2.1);
  
\end{tikzpicture}
\caption{Illustration of the fine grid, coarse grid $K_i$ and oversampling domain $K_{i}^{+}$.}
\label{fig:grid2}
\end{figure}

\section{Multiscale space construction}
\label{section:Multiscale space construction}
In this section, we first briefly describe the construction of multiscale spaces based on CEM-GMsFEM \cite{reference9,reference14,reference15}. We then present the construction of the spaces $V_{H,1}$ and $V_{H,2}$ based on the non-local multi-continuum (NLMC) method \cite{reference11,reference33}. 
\subsection{CEM-GMsFEM}
Denote by $\{K_{i}\}$ the set of coarse blocks in $\mathcal{T}_{H}$, and denote $V=H_{0}^{1}(\Omega)$. For each $K_{i}\in\mathcal{T}_{H}$, we have to build a collection of auxiliary based in $V(K_{i})$ where $V(K_i)$ is the restriction of $V$ on $K_i$. We solve the following eigenvalue problem 
\begin{align}\label{eqnBa}
\int_{K_{i}}\kappa\nabla\psi_{j}^{(i)}\cdot\nabla v=\lambda_{j}^{(i)}s_{i}\left(\psi_{j}^{(i)},v\right),\ \forall v\in V(K_{i}),
\end{align}
where
\begin{align*}
s_{i}(u,v)=\int_{K_{i}}\widetilde{\kappa}uv,\ \widetilde{\kappa}=\kappa\sum_{i}\mid\nabla\chi_{i}\mid^{2}\ \text{or}\ \widetilde{\kappa}=\kappa H^{-2},
\end{align*}
and $\{\chi_{i}\}$ is the partition of unity corresponding to an overlapping partition of the domain. Then we collect the first $J_{i}$ eigenfunctions corresponding to the first $J_{i}$ smallest eigenvalues to form the auxiliary spaces
\begin{align*}
V_{\text{aux}}^{(i)}:=\text{span}\{\psi_{j}^{(i)}:1\leq j\leq J_{i}\}.
\end{align*}
Define a projection operator $\Pi: L^{2}(\Omega)\rightarrow V_{\text{aux}}\subset L^{2}(\Omega)$
\begin{align*}
s\left(\Pi u,v\right)=s\left(u,v\right),\ \forall v\in V_{\text{aux}}:=\sum_{i=1}^{N_{e}}V_{\text{aux}}^{(i)},
\end{align*}
where $s(u,v)=\sum_{i=1}^{N_{e}}s_{i}(u|_{K_{i}},v|_{K_{i}})$ and $N_{e}$ is the number of coarse elements. Define $K_{i}^{+}$ to be an oversampling domain obtained by enlarging $K_{i}$ by a suitable number of coarse grid layers, see Figure \ref{fig:grid2} for an illustration. For each $\psi_{j}^{(i)}$, we search for a local basis function $\phi_{j}^{(i)}\in V(K_{i}^{+})$ such that for some $\mu_{j}^{(i)}\in V(K_{i}^{+})$
\begin{align}\label{eqnBb}
\begin{aligned}
a\left(\phi_{j}^{(i)},v\right)+s\left(\mu_{j}^{(i)},v\right)&=0,\ \forall v\in V(K_{i}^{+}),\\
s\left(\phi_{j}^{(i)},\nu\right)&=s\left(\psi_{j}^{(i)},\nu\right),\ \forall\nu\in V_{\text{aux}}(K_{i}^{+}).
\end{aligned}
\end{align}
Then we define the space $V_{cem}$ as
\begin{align}\label{eqnBc}
V_{cem}:=\text{span}\{\phi_{j}^{(i)}: 1\leq i\leq N_{e}, 1\leq j\leq J_{i}\}.
\end{align}
Thus we can choose two subspaces $V_{H,1}$ and $V_{H,2}$ based on the multiscale space $V_{{cem}}$ \cite{reference5,reference6}. %in $\widetilde{V}$ 

\subsection{Construction of two multiscale subspaces}

For channelized media, the construction of multiscale subspaces can be simplified. That is, denote by the computational domain $\Omega = \Omega_m \bigoplus_{l=1}^s d_l\Omega_{f,l}$, where $m$ and $f$ denote the matrix region and fracture region. In the fracture regions $\Omega_{f,l}$, the scalar $d_{l}$ and $s$ denote the aperture and the number of the discrete fracture networks, respectively. Since the value of the permeabilities in the matrix and fracture regions can differ in magnitudes, thus we can construct constraint energy minimizing basis functions via NLMC, such that the obtained basis functions can automatically separate two continua. Specifically, for a given coarse block $K_{i}$, we use constants for each individual fracture network and then a constant for the matrix to form a auxiliary space. That is to say, for any coarse block $K_{i}$, we write $K_{i}=K_{i,f}\cup K_{i,m}$, where $K_{i,f}:=\{f_{j}^{(i)}, j=1,\cdots,m_{i}\}$ is the high-contrast channelized region, $m_{i}$ is the number of non-connected fractures in $K_{i}$, $K_{i,m}$ is its complement in $K_{i}$. Then we define two auxiliary spaces
\begin{align}
\label{eqnBd}&V_{\text{aux},1}^{(i)}:=\text{span}\{\phi_{\text{aux},k}^{(i)}\mid\phi_{\text{aux},k}^{(i)}=0\ \text{in}\ K_{i,m},\phi_{\text{aux},k}^{(i)}=\delta_{jk}\ \text{in}\ f_{j}^{(i)}, k=1,\cdots,m_{i}\},\\
\label{eqnBe}&V_{\text{aux},2}^{(i)}:=\text{span}\{\phi_{\text{aux},0}^{(i)}\mid\phi_{\text{aux},0}^{(i)}=1\ \text{in}\ K_{i,m},\phi_{\text{aux},0}^{(i)}=0\ \text{in}\ K_{i,f}\}.
\end{align}
\indent Then the NLMC basis functions are obtained by finding $\psi_{m}^{(i)}\in V_{0}(K_{i}^{+})$ and $\mu_{0}^{(j)},\mu_{n}^{(j)}\in\mathbb{R}$ from the following localized constraint energy minimizing problem
\begin{align}\label{eqnBf}
\left\{\begin{aligned}
&a\left(\psi_{m}^{(i)},v\right)+\sum_{K_{j}\subset K_{i}^{+}}\left(\mu_{0}^{(j)}\int_{K_{j,m}}v+\sum_{1\leq n\leq m_{j}}\mu_{n}^{(j)}\int_{f_{n}^{(j)}}v\right)=0,\ \forall v\in V_{0}(K_{i}^{+}),\\
&\int_{K_{j,m}}\psi_{m}^{(i)}=\delta_{ij}\delta_{m0},\ \forall K_{j}\subset K_{i}^{+}, 0\leq m\leq m_i,\\
&\int_{f_{n}^{(i)}}\psi_{m}^{(i)}=\delta_{ij}\delta_{mn},\ \forall f_{n}^{(j)}\in\mathcal{F}_{j},\ \forall K_{j}\subset K_{i}^{+}.
\end{aligned}\right.
\end{align}

The above problems are posed in infinite-dimensional spaces, but for numerical computations, we solve the discretized system on the fine grid using standard finite elements to obtain the solutions and use them as our basis $\{\psi_{m}^{(i)},0\leq m\leq m_{i}, 1\leq i\leq N_{e}\}$.

\indent Denote the average of all NLMC basis by
\begin{align}\label{eqnBg}
\overline{\psi}:=\frac{1}{L}\sum_{i=1}^{N_{e}}\sum_{m=0}^{m_{i}}\psi_{m}^{(i)},\quad L=\sum_{i=1}^{N}m_{i}.
\end{align}
We then let $\widetilde{\psi}_{m}^{(i)}=\psi_{m}^{(i)}-\frac{s(\psi_{m}^{(i)},\overline{\psi})}{s(\overline{\psi},\overline{\psi})}\overline{\psi},0\leq m\leq m_{i}, 1\leq i\leq N_{e}$. In order to simplify the notations, we omit the double script in $\widetilde{\psi}_{m}^{(i)}$ and denote the set of bases by $\{\widetilde{\psi}_{k},k=1,\cdots,L\}$. Thus, we define the space $V_{H,1}$ as
\begin{align*}
V_{H,1}=\text{span}\{\widetilde{\psi}_{k},1\leq k\leq L-1\}.
\end{align*}
The basis functions corresponding to the matrix and the basis $\overline{\psi}$ will be included in the second subspace $V_{H,2}$, that is,
\begin{align*}
V_{H,2}=\text{span}\{\overline{\psi},\psi_{0}^{(i)},1\leq i\leq N_{c}\}.
\end{align*}
\indent By this construction, $V_{H,1}$ contains a basis representing the high-contrast fractures only, and $V_{H,2}$ includes a basis representing the background matrix and the constant basis. With this choice of $V_{H,1}$ and $V_{H,2}$, one can show that the splitting scheme \eqref{eqnAd}-\eqref{eqnAe} is stable with the condition $\Delta t \sup_{v\in V_{H,2}}\cfrac{\|v\|_a^2}{\|v\|^2} \leq 1-\gamma^2$ \cite{reference14, reference32}. %Noting that we take away the last basis in $V_{H,1}$ to remove linear dependency between the two space.  

\indent Next, we will introduce some notations. Let dim$(V_{H,1})=d_{1}$, dim$(V_{H,2})=d_{2}$, dim$(V_{H})=D$, and let $\Psi_{1}\in\mathbb{R}^{D\times d_{1}}$ and $\Psi_{2}\in\mathbb{R}^{D\times d_{2}}$ be the matrices whose columns are the bases of $V_{H,1}$ and $V_{H,2}$, respectively. Denote $M_{f}$ and $A_{f}$ be the fine scale mass matrix and stiffness matrix, define the following coarse scale matrices
\begin{align}\label{eqnBh}
\begin{aligned}
M_{11}&=\Psi_{1}^{T}M_{f}\Psi_{1},\ A_{11}=\Psi_{1}^{T}A_{f}\Psi_{1},\\
M_{22}&=\Psi_{2}^{T}M_{f}\Psi_{2},\ A_{22}=\Psi_{2}^{T}A_{f}\Psi_{2},\\
M_{12}&=\Psi_{1}^{T}M_{f}\Psi_{2},\ A_{12}=\Psi_{1}^{T}A_{f}\Psi_{2},\\
F_{1}^{n}&=\Psi_{1}^{T}f^{n},\quad\quad~ F_{2}^{n}=\Psi_{2}^{T}f^{n}.
\end{aligned}
\end{align}

\section{Parallel in time for partially explicit temporal splitting scheme}
\label{section:Parallel in time for partially explicit temporal splitting scheme}
In this section, we describe in detail the parareal algorithm for the partially explicit temporal splitting scheme. First, we present the basis flow of the WR technique and the diagonalization in section \ref{section: WR technique and diagonalization technique}. Then we introduce the WR method via the diagonalization technique based on the all-at-once system for the partially explicit scheme in subsection \ref{section:WR method via diagonalization for all-at-once system} and give the convergence of the method. Then we propose our main algorithm, the parareal all-at-once partially explicit temporal splitting algorithm, in section \ref{section:Parareal all-at-once partially explicit temporal splitting algorithm}, where we adopt the WR method via diagonalization technique as the fine solver. Finally, we carry out the error analysis for the proposed algorithm. 

\subsection{WR technique and diagonalization technique}
\label{section: WR technique and diagonalization technique}
We first present the WR technique \cite{reference39,reference40,reference41} and the diagonalization \cite{reference40,reference41} by considering the initial value problem \eqref{eqnAf}. The WR technique can be summarized as follows. At first, noticing that $u(0)=\alpha u(T)-\alpha u(T)+u_0$ holds for all $\alpha\in\mathbb{R}$, then we can construct the following WR iterations
\begin{align}\label{WR1}
\left\{\begin{aligned}
&\frac{du^j}{dt}=F(t,u^j),\ t\in(0,T),\\
&u^j(0)=\alpha u^j(T)-\alpha u^{j-1}(T)+u_0,
\end{aligned}\right.
\end{align}
where $j\geq1$ is the iteration index and $u^0(t)$ denotes the initial guess. Then for each iteration of \eqref{WR1}, it is a differential equation with periodic-like condition. Thus the diagonalization technique can be adopted and carry out parallel computation. To illustrate it in detail, we consider $F(t,u)=-Au(t)+\widetilde{f}(t)$ and utilize the Backward Euler scheme to discretize the differential equation \eqref{WR1}, then we obtain the following 
\begin{align}\label{WR2}
\left\{\begin{aligned}
&\frac{u_{n+1}^j-u_n^j}{\Delta t}=-Au_{n+1}^j+\widetilde{f}_{n+1},\\
&u_0^j=\alpha u_N^j+R^{j-1}\ \text{with}\ R^{j-1}:=-\alpha u_N^{j-1}+u_0,
\end{aligned}\right.
\end{align}
where $n=1,2,\cdots,N:=T/\Delta t$. Next, one can rewrite \eqref{WR2} as the following all-at-once system
\begin{align}\label{WR3}
(B\otimes I_N+I_N\otimes A)U^j=\widetilde{F},
\end{align}
where $U^j=\left(u_1^j,u_2^j,\cdots,u_N^j\right)^{T}$, $\widetilde{F} = [\frac{1}{\Delta t}(u_0-\alpha u_N^{j-1})+\widetilde{f}_1, \widetilde{f}_2, \cdots, \widetilde{f}_n]^T$, $B\in\mathbb{R}^{N\times N}$ is a periodic-like matrix with parameter $\alpha\in(0,1)$ as follows
\begin{align}\label{WR4}
B=\frac{1}{\Delta t}\left(
  \begin{array}{cccc}
    1 & ~ & ~ & -\alpha \\
    -1 & 1 & ~ & ~ \\
    ~ & \ddots & \ddots & ~ \\
    ~ & ~ & -1 & 1 \\
  \end{array}
\right).
\end{align}
The $\alpha$-circulant matrix $B$ can then be diagonalized, i.e.
\begin{align}\label{WR5}
B=SDS^{-1},\ D=\text{diag}(d_1,d_2,\cdots,d_N),    
\end{align}
then we can factorize the coefficient matrix in \eqref{WR3} as follows
\begin{align}\label{WR6}
B\otimes I_N+I_N\otimes A=(S\otimes I_N)(D\otimes I_N+I_N\otimes A)(S^{-1}\otimes I_N).
\end{align}
Then we can solve \eqref{WR3} in the $j-$th iteration in the following three steps.
\begin{align}\label{WR7}
\begin{aligned}
&(a) ~\left(S\otimes I_{N}\right)P=F^{j},\\
&(b) ~\left(D\otimes I_N+I_{N}\otimes A\right)Q=P,\\
&(c) ~\left(S^{-1}\otimes I_{N}\right)U^{j} = Q.
\end{aligned}
\end{align}
Note that the matrix $S$ can be further decomposed into $S=\Lambda V$, where $\Lambda=\text{diag}\{1,\alpha^{-\frac{1}{N}},\cdots,\alpha^{-\frac{N-1}{N}}\}$, $V$ is the discrete Fourier matrix. Thus, the Fast Fourier Transform (FFT) can be employed to speed up the implementation of \eqref{WR7}. In addition, the second step of \eqref{WR7} is to solve $N$ independent equations, thus it can be done in parallel.

\subsection{WR method via diagonalization for all-at-once system}
\label{section:WR method via diagonalization for all-at-once system}
Let us look back at the splitting scheme \eqref{eqnAd}-\eqref{eqnAe}, and solve it using the iterative all-at-once method. We write
\begin{align}
\label{WR8}&\left(\frac{u_{n+1}^{j}-u_{n}^{j}}{\Delta t},v_{1}\right)+a\left(u_{n+1}^{j},v_{1}\right)=\left(f^{n+1},v_{1}\right)-\left(\frac{w_{n}^{j-1}-w_{n-1}^{j-1}}{\Delta t},v_{1}\right)-a\left(w_{n}^{j-1},v_{1}\right),\\
\label{WR9}&\left(\frac{w_{n+1}^{j}-w_{n}^{j}}{\Delta t},v_{2}\right)+\left(\frac{u_{n}^{j}-u_{n-1}^{j}}{\Delta t},v_{2}\right)+a\left(u_{n+1}^{j}+w_{n}^{j},v_{2}\right)=\left(f^{n+1},v_{2}\right),
\end{align}
where $v_{1}\in V_{H,1}, v_{2}\in V_{H,2}$ and $j$ will be the iteration index. Let $U^{j}=(u_{1}^{j},u_{2}^{j},\cdots,u_{N}^{j})$ and $W^{j}=(w_{1}^{j},w_{2}^{j},\cdots,w_{N}^{j})$, therefore \eqref{WR8}-\eqref{WR9} is a direct discretization of the following semi-discretization scheme: $\forall v_{1}\in V_{H,1}$ and $\forall v_{2}\in V_{H,2}$
\begin{align}
\label{WR10}&\left(\partial_{t}U^{j}(t),v_{1}\right)+a\left(U^{j}(t),v_{1}\right)=\left(f,v_1\right)-\left(\partial_{t}W^{j-1}(t),v_{1}\right)-a\left(W^{j-1}(t),v_{1}\right),\\
\label{WR11}&\left(\partial_{t}W^{j}(t),v_{2}\right)+a\left(W^{j}(t),v_{2}\right)=\left(f,v_2\right)-\left(\partial_{t}U^{j}(t),v_{2}\right)-a\left(U^{j}(t),v_{2}\right).
\end{align}
Then we rewrite the \eqref{WR8} into the following all-at-once system
\begin{align}\label{WR12}
\left(B\otimes M_{11}+I_{t}\otimes A_{11}\right)U^{j}=F
\end{align}
with initial condition 
\begin{align*}
U^{j}(0)= u_{0}+\alpha\left(u_N^{j}-u_N^{j-1}\right),
\end{align*}
where $B$ is a periodic-like matrix defined in \eqref{WR4}, $F=[\widetilde{f}^1, \cdots, \widetilde{f}^N]^T$ and
$$\widetilde{f}^{1}=F_{1}^{1}-\frac{1}{\Delta t}M_{12}w_{0}^{j-1}-A_{12}w_{0}^{j-1}+\frac{1}{\Delta t}M_{11}\left(u_{0}-\alpha u_{N}^{j-1}\right),$$
$$\widetilde{f}^{s}=F_{1}^{s}-\frac{1}{\Delta t}M_{12}\left(w_{s}^{j-1}-w_{s-1}^{j-1}\right)-A_{12}w_{s}^{j-1},$$ for $s=2,3,\cdots,N$. Similar to the section \ref{section: WR technique and diagonalization technique}, we assume that the matrix $B$ is diagonalizable and \eqref{WR5} holds, and factorize the coefficient matrix in \eqref{WR12} as follows
\begin{align}\label{WR13}
B\otimes M_{11}+I_{t}\otimes A_{11}=\left(S\otimes I_{M}\right)\left(D\otimes M_{11}+I_{t}\otimes A_{11}\right)\left(S^{-1}\otimes I_{M}\right),
\end{align}
where $I_{M}$ is an identity matrix. Then we can solve the system \eqref{WR12} in three steps like \eqref{WR7}, and adopt the FFT to speed up the computation.

The problem \eqref{WR8}-\eqref{WR9} will be solved with the following steps:
\begin{enumerate}
    \item Solve \eqref{WR12} by corresponding three steps \eqref{WR7}, obtain $U^{j}=(u_{1}^{j},\cdots,u_{N}^{j})$, which is the part of the solution in $V_{H,1}$ at all time steps.
    \item Plug the solution $U$ in \eqref{WR9}, solve for $W$ in a sequential manner.
    \item Iterate the above process until converge.
\end{enumerate}

We consider
\begin{align*}
&\left(\partial_{t}U(t),v_{1}\right)+a\left(U(t),v_{1}\right)=-\left(\partial_{t}W(t),v_{1}\right)-a\left(W(t),v_{1}\right),\ \forall v_{1}\in V_{H,1},\\
&\left(\partial_{t}W(t),v_{2}\right)+a\left(W(t),v_{2}\right)=-\left(\partial_{t}U(t),v_{2}\right)-a\left(U(t),v_{2}\right),\ \forall v_{2}\in V_{H,2}
\end{align*}
The following theorem gives the convergence result for the WR method at each subinterval $(t_n, t_{n+1})$. We remark that we will adopt the WR method as the fine propagator at each subinterval $(t_n, t_{n+1})$ within the parareal framework. %It is trivial to extend the following results on the whole temporal domain $(0,T]$. 
\begin{theorem}
Let $U^{j}(t)$ and $W^{j}(t)$ are the solution of \eqref{WR10}-\eqref{WR11} at $j$-th iteration respectively. Then it holds for $j=1,2,\cdots$ and $n=0,1,2,\cdots,N-1$ that
\begin{align}\label{WR14}
\begin{aligned}
&\sup_{t\in(t_{n},t_{n+1})}\parallel U^{j}(t)-U(t)\parallel+\sup_{t\in(t_{n},t_{n+1})}\parallel W^{j}(t)-W(t)\parallel\\
&\leq C\Delta t\gamma^{2(j-1)}\sup_{t\in(t_{n},t_{n+1})}\parallel A_{22}^{-1}\partial_{t}(W^{0}(t)-W(t))+(W^{0}(t)-W(t))\parallel_{\infty},
\end{aligned}
\end{align}
where $\gamma$ is defined by 
\begin{align}\label{WR15}
\gamma=\frac{|(v_{1},v_{2})|}{\parallel v_{1}\parallel\parallel v_{2}\parallel},\ \forall v_{1}\in V_{H,1},\ \forall v_{2}\in V_{H,2}.
\end{align}
\end{theorem}
\textbf{Proof:}
\indent We define $P_{1}: V_{H,2}\rightarrow V_{H,1}$, $\Pi_{1}: V_{H,2}\rightarrow V_{H,1}$ and $P_{2}: V_{H,1}\rightarrow V_{H,2}$, $\Pi_{2}: V_{H,1}\rightarrow V_{H,2}$ such that
\begin{align*}
\left(P_{1}v_{2},v_{1}\right)=&\left(v_{2},v_{1}\right),\ \forall v_{1}\in V_{H,1},\ \left(P_{2}v_{1},v_{2}\right)=\left(v_{1},v_{2}\right),\ \forall v_{2}\in V_{H,2},\\
a\left(\Pi_{1}v_{2},v_{1}\right)=&a\left(v_{2},v_{1}\right),\ \forall v_{1}\in V_{H,1},\ a\left(\Pi_{2}v_{1},v_{2}\right)=a\left(v_{1},v_{2}\right),\ \forall v_{2}\in V_{H,2},
\end{align*}
and define the errors
\begin{align*}
e_{u}^{j}(t)=U^{j}(t)-U(t),\ e_{w}^{j}(t)=W^{j}(t)-W(t).
\end{align*}
Then for $\forall v_{1}\in V_{H,1}$ and $\forall v_{2}\in V_{H,2}$, we have
\begin{align}
\label{WR16}&\left(\partial_{t}e_{u}^{j}(t),v_{1}\right)+a\left(e_{u}^{j}(t),v_{1}\right)=-\left(\partial_{t}e_{w}^{j-1}(t),P_{2}v_{1}\right)-a\left(e_{w}^{j-1}(t),\Pi_{2}v_{1}\right),\\
\label{WR17}&\left(\partial_{t}e_{w}^{j}(t),v_{2}\right)+a\left(e_{w}^{j}(t),v_{2}\right)=-\left(\partial_{t}e_{u}^{j}(t),P_{1}v_{2}\right)-a\left(e_{u}^{j}(t),\Pi_{1}v_{2}\right).
\end{align}
with initial condition $e_{u}^{j}(0)=e_{w}^{j}(0)=0$. Then from \eqref{WR16}-\eqref{WR17} we obtain
\begin{align}\label{WR18}
&\left(\partial_{t}e_{w}^{j}(t),v_{2}\right)+a\left(e_{w}^{j}(t),v_{2}\right)+a\left(e_{u}^{j}(t),\Pi_{1}v_{2}\right)\nonumber\\
=&a\left(e_{u}^{j}(t),P_{1}v_{2}\right)+\left(\partial_{t}e_{w}^{j-1}(t),P_{2}P_{1}v_{2}\right)+a\left(e_{w}^{j-1}(t),\Pi_{2}P_{1}v_{2}\right).
\end{align}
From \eqref{WR17}, a directly calculation yields the following recurrence relation
\begin{align*}
&-\left(\partial_{t}e_{w}^{j-1}(t),v_{2}\right)\\
=&a\left(e_{w}^{j-1}(t),v_{2}\right)+\left(\partial_{t}e_{u}^{j-1}(t),P_{1}v_{2}\right)+a\left(e_{u}^{j-1}(t),\Pi_{1}v_{2}\right)\\
=&a\left(e_{w}^{j-1}(t),v_{2}\right)-a\left(e_{u}^{j-1}(t),P_{1}v_{2}\right)-\left(\partial_{t}e_{w}^{j-2}(t),P_{2}P_{1}v_{2}\right)\\
&-a\left(e_{w}^{j-2}(t),\Pi_{2}P_{1}v_{2}\right)+a\left(e_{u}^{j-1}(t),\Pi_{1}v_{2}\right)\\
=&a\left(e_{w}^{j-1}(t),v_{2}\right)-a\left(e_{w}^{j-2}(t),\Pi_{2}P_{1}v_{2}\right)+a\left(e_{u}^{j-1}(t),\Pi_{1}v_{2}-P_{1}v_{2}\right)\\
&-\left(\partial_{t}e_{w}^{j-2}(t),P_{2}P_{1}v_{2}\right)\\
=&a\left(e_{w}^{j-1}(t),v_{2}\right)+a\left(e_{u}^{j-1}(t),\Pi_{1}v_{2}-P_{1}v_{2}\right)+a\left(e_{u}^{j-1}(t),P_{2}P_{1}v_{2}-\Pi_{2}P_{1}v_{2}\right)\\
&+a\left(e_{u}^{j-2}(t),\Pi_{1}P_{2}P_{1}v_{2}-P_{1}P_{2}P_{1}v_{2}\right)-a\left(e_{w}^{j-3}(t),\Pi_{2}P_{1}P_{2}P_{1}v_{2}\right)\\
&-\left(\partial_{t}e_{w}^{j-3}(t),P_{2}P_{1}P_{2}P_{1}v_{2}\right).
\end{align*}
Denote by $P_{21}=P_{2}P_{1}$ and $P_{12}=P_{1}P_{2}$, then the following recurrence relation can be written as
\begin{align*}
&-\left(\partial_{t}e_{w}^{j-1}(t),v_{2}\right)=a\left(e_{w}^{j-1}(t),v_{2}\right)+a\left(e_{u}^{j-1}(t),\Pi_{1}v_{2}-P_{1}v_{2}\right)\\
+&\sum_{i=0}\left(a(e_{w}^{j-2-i}(t),P_{21}^{i+1}v_{2}-\Pi_{2}P_{1}P_{21}^{i}v_{2})+a(e_{u}^{j-2-i}(t),\Pi_{1}P_{21}^{i+1}v_{2}-P_{1}P_{21}^{i+1}v_{2})\right)\\
-&a\left(e_{w}^{0}(t),\Pi_{2}P_{1}P_{21}^{j-2}v_{2}\right)-\left(\partial_{t}e_{w}^{0}(t),P_{21}^{j-1}v_{2}\right).
\end{align*}
Substitute the above equation into \eqref{WR16} and \eqref{WR18} respectively, we have
\begin{align}\label{WR19}
\begin{aligned}
&\left(\partial_{t}e_{u}^{j}(t),v_{1}\right)+a\left(e_{u}^{j}(t),v_{1}\right)\\
=&\sum_{i=0}\left(a(e_{w}^{j-1-i}(t),(P_{2}-\Pi_{2})P_{12}^{i}v_{1})+a(e_{u}^{j-1-i}(t),(\Pi_{1}-P_{1})P_{2}P_{12}^{i}v_{1})\right)\\
&-a\left(e_{w}^{0}(t),\Pi_{2}P_{12}^{j-1}v_{1}\right)-\left(\partial_{t}e_{w}^{0}(t),P_{2}P_{12}^{j-1}v_{1}\right)
\end{aligned}
\end{align}
and
\begin{align}\label{WR20}
\begin{aligned}
&\left(\partial_{t}e_{w}^{j}(t),v_{2}\right)+a\left(e_{w}^{j}(t),v_{2}\right)+a\left(e_{u}^{j}(t),\Pi_{1}v_{2}\right)\\
=&\sum_{i=0}\left(a(e_{w}^{j-1-i}(t),(\Pi_{2}-P_{2})P_{1}P_{21}^{i}v_{2})+a(e_{u}^{j-1-i}(t),(P_{1}-\Pi_{1})P_{21}^{i+1}v_{2})\right)\\
&+a\left(e_{w}^{0}(t),\Pi_{2}P_{1}P_{21}^{j-1}v_{2}\right)+\left(\partial_{t}e_{w}^{0}(t),P_{21}^{j}v_{2}\right).
\end{aligned}
\end{align}
\indent We then write the above operator into the matrix form
\begin{align}\label{WR21}
\begin{aligned}
&\partial_{t}e_{u}^{j}(t)+A_{11}e_{u}^{j}(t)\\
=&\sum_{i=0}\left(P_{12}^{i*}(P_{2}^{*}-\Pi_{1})A_{22}e_{w}^{j-1-i}(t)+P_{12}^{i*}P_{2}^{*}(\Pi_{2}-P_{1}^{*})A_{11}e_{u}^{j-1-i}(t)\right)\\
&-\left(P_{12}^{j-1}\right)^{*}\Pi_{1}A_{22}e_{w}^{0}(t)-\left(P_{12}^{j-1}\right)^{*}P_{2}^{*}\partial_{t}e_{w}^{0}(t)
\end{aligned}
\end{align}
and
\begin{align}\label{WR22}
\begin{aligned}
&\partial_{t}e_{w}^{j}(t)+A_{22}e_{w}^{j}(t)+\left(\Pi_{2}-P_{1}^{*}\right)A_{11}e_{u}^{j}(t)\\
=&\sum_{i=0}\left(P_{21}^{i*}P_{1}^{*}(P_{2}^{*}-\Pi_{1})A_{22}e_{w}^{j-1-i}(t)+(P_{21}^{i+1})^{*}(\Pi_{2}-P_{1}^{*})A_{11}e_{u}^{j-1-i}(t)\right)\\
&+\left(P_{21}^{j-1}\right)^{*}P_{1}^{*}\Pi_{1}A_{22}e_{w}^{0}(t)+\left(P_{21}^{j}\right)^{*}\partial_{t}e_{w}^{0}(t).
\end{aligned}
\end{align}
Denote by
\begin{align*}
f_{w}^{j}(t)=&\sum_{i=0}P_{12}^{i*}\left(P_{2}^{*}-\Pi_{1}\right)A_{22}e_{w}^{j-1-i}(t),\\
f_{u}^{j}(t)=&\sum_{i=0}P_{12}^{i*}P_{2}^{*}\left(\Pi_{2}-P_{1}^{*}\right)A_{11}e_{u}^{j-1-i}(t),
\end{align*}
then with the above two notations, we can rewrite \eqref{WR21} and \eqref{WR22} as follows
\begin{align}\label{WR23}
\partial_{t}e_{u}^{j}(t)+A_{11}e_{u}^{j}(t)=f_{w}^{j}(t)+f_{u}^{j}(t)-\left(P_{12}^{j-1}\right)^{*}\Pi_{1}A_{22}e_{w}^{0}(t)-\left(P_{12}^{j-1}\right)^{*}P_{2}^{*}\partial_{t}e_{w}^{0}(t)
\end{align}
and 
\begin{align}\label{WR24}
\begin{aligned}
&\partial_{t}e_{w}^{j}(t)+A_{22}e_{w}^{j}(t)+(\Pi_{2}-P_{1}^{*})A_{11}e_{u}^{j}(t)\\
=&P_{1}^{*}f_{w}^{j}(t)+P_{1}^{*}f_{u}^{j}(t)-P_{1}^{*}(P_{12}^{j-1})^{*}\Pi_{1}A_{22}e_{w}^{0}(t)-P_{1}^{*}(P_{12}^{j-1})^{*}P_{2}^{*}\partial_{t}e_{w}^{0}(t).
\end{aligned}
\end{align}
We consider $\parallel A_{11}^{-1}P_{2}^{*}A_{22}\parallel\leq C$ and $\parallel A_{22}^{-1}P_{1}^{*}A_{11}\parallel\leq C$, where $C$ is independent of contrast, and let $\lambda_{p}$ and $\phi_{p}$ be the eigenvalues and eigenfunctions such that
\begin{align*}
P_{12}^{*}\phi_{p}=\lambda_{p}\phi_{p},\ p=1,2,\cdots.
\end{align*}
Thus, for $\widetilde{\phi}_{p}=A_{11}^{-1}\phi_{p}$, we have $A_{11}^{-1}P_{12}^{*}A_{11}\widetilde{\phi}_{p}=\lambda_{i}\widetilde{\phi}_{p}$. Let $\widetilde{\Phi}$ and $D$ be the matrix $\widetilde{\Phi}=[\widetilde{\phi}_{1},\widetilde{\phi}_{2},\cdots]$ and $D=\text{diag}\{\lambda_{1},\lambda_{2},\cdots\}$, then 
\begin{align*}
f_{w}^{j}(t)=&\sum_{i=0}P_{12}^{i*}\left(P_{2}^{*}-\Pi_{1}\right)A_{22}e_{w}^{j-1-i}(t)\\
=&\sum_{i=0}P_{12}^{i*}A_{11}\left(A_{11}^{-1}P_{2}^{*}A_{22}-\Pi_{1}\right)e_{w}^{j-1-i}(t)\\
=&\sum_{i=0}P_{12}^{i*}A_{11}\widetilde{\Phi}\widetilde{\Phi}^{-1}\left(A_{11}^{-1}P_{2}^{*}A_{22}-\Pi_{1}\right)e_{w}^{j-1-i}(t)\\
=&\sum_{i=0}A_{11}\widetilde{\Phi}D^{i}\widetilde{\Phi}^{-1}\left(A_{11}^{-1}P_{2}^{*}A_{22}-\Pi_{1}\right)e_{w}^{j-1-i}(t)
\end{align*}
Since $\lambda_{p}\leq\gamma^{2}, p=1,2,\cdots$, where $\gamma=\frac{|(v_{1},v_{2})|}{\parallel v_{1}\parallel\parallel v_{2}\parallel}$, then we have
\begin{align}\label{WR25}
\parallel A_{11}^{-1}f_{w}^{j}(t)\parallel=&\parallel\sum_{i=0}\widetilde{\Phi}D^{i}\widetilde{\Phi}^{-1}\left(A_{11}^{-1}P_{2}^{*}A_{22}-\Pi_{1}\right)e_{w}^{j-1-i}(t)\parallel\nonumber\\
\leq&C\sum_{i=0}\max_{p}\{\lambda_{p}\}^{i}\parallel e_{w}^{k-1-i}\parallel\nonumber\\
\leq&C\frac{1}{1-\gamma^{2}}\max_{0\leq i\leq j-1}\parallel e_{w}^{i}(t)\parallel.
\end{align}
Similarly we have
\begin{align}
\label{WR26}\parallel A_{11}^{-1}f_{u}^{j}(t)\parallel\leq C\frac{1}{1-\gamma^{2}}\max_{0\leq i\leq j-1}\parallel e_{u}^{i}(t)\parallel,\\
\label{WR27}\parallel A_{22}^{-1}P_{1}^{*}f_{w}^{j}(t)\parallel\leq C\frac{1}{1-\gamma^{2}}\max_{0\leq i\leq j-1}\parallel e_{w}^{i}(t)\parallel,\\
\label{WR28}\parallel A_{22}^{-1}P_{1}^{*}f_{u}^{j}(t)\parallel\leq C\frac{1}{1-\gamma^{2}}\max_{0\leq i\leq j-1}\parallel e_{w}^{i}(t)\parallel.
\end{align}
Then we turn to estimate the last two terms of \eqref{WR21} and \eqref{WR22},we obtain
\begin{align}
\label{WR29}&\parallel A_{11}^{-1}\left((P_{12}^{j-1})^{*}\Pi_{1}A_{22}e_{w}^{0}(t)+(P_{12}^{j-1})^{*}P_{2}^{*}\partial_{t}e_{w}^{0}(t)\right)\parallel\nonumber\\
\leq&C\gamma^{2(k-1)}\parallel A_{11}^{-1}(P_{2}^{*}\partial_{t}e_{w}^{0}(t)+\Pi_{1}A_{22}e_{w}^{0}(t))\parallel,\\
\label{WR30}&\parallel A_{22}^{-1}P_{1}^{*}\left((P_{12}^{j-1})^{*}\Pi_{1}A_{22}e_{w}^{0}(t)+(P_{12}^{j-1})^{*}P_{2}^{*}\partial_{t}e_{w}^{0}(t)\right)\parallel\nonumber\\
\leq&C\gamma^{2(k-1)}\parallel A_{22}^{-1}P_{1}^{*}\left(P_{2}^{*}\partial_{t}e_{w}^{0}(t)+\Pi_{1}A_{22}e_{w}^{0}(t)\right)\parallel.
\end{align}
Now we consider the following matrices
\begin{align*}
\widetilde{A}=\left(
  \begin{array}{cc}
    A_{11} & O \\
    \left(\Pi_{2}-P_{1}^{*}\right)A_{11} & A_{22} \\
  \end{array}
\right),\ 
A=\left(
  \begin{array}{cc}
    A_{11} & O \\
    O & A_{22} \\
  \end{array}
\right)
\end{align*}
then we have
\begin{align*}
\widetilde{A}^{-1}=\left(
  \begin{array}{cc}
    A_{11}^{-1} & O \\
    -A_{22}^{-1}\left(\Pi_{2}-P_{1}^{*}\right) & A_{22}^{-1} \\
  \end{array}
\right),
\widetilde{A}^{-1}A=\left(
  \begin{array}{cc}
    I & O \\
    -A_{22}^{-1}\left(\Pi_{2}-P_{1}^{*}\right)A_{11} & I \\
  \end{array}
\right)
\end{align*}
Thus, we can easily get $\parallel\widetilde{A}^{-1}A\parallel\leq C$ and 
\begin{align*}
\parallel\int_{0}^{t}e^{-(t-s)\widetilde{A}}A\parallel=\parallel\left(I-e^{-t\widetilde{A}}\right)\widetilde{A}^{-1}A\parallel\leq\parallel I-e^{-t\widetilde{A}}\parallel\parallel \widetilde{A}^{-1}A\parallel\leq Ct.
\end{align*}
Therefore, given any $\Delta t>0$, $t_{n}=n\Delta t$, and with the help of the \eqref{WR25}-\eqref{WR30}, we have the following estimates
\begin{align}\label{WR31}
\sup_{t\in(t_{n},t_{n+1})}&\parallel e_{u}^{j}(t)\parallel\leq\parallel e_{u}^{j}(t_{n})\parallel+C\Delta t\gamma^{2(j-1)}\parallel A_{11}^{-1}(P_{2}^{*}\partial_{t}e_{w}^{0}(t)+\Pi_{1}A_{22}e_{w}^{0}(t))\parallel\nonumber\\
+&C\Delta t\frac{1}{1-\gamma^{2}}\left(\max_{0\leq i\leq j-1}\sup_{t\in(t_{n},t_{n+1})}\parallel e_{w}^{i}(t)\parallel+\max_{0\leq i\leq j-1}\sup_{t\in(t_{n},t_{n+1})}\parallel e_{u}^{i}(t)\parallel\right),
\end{align}
\begin{align}\label{WR32}
\sup_{t\in(t_{n},t_{n+1})}&\parallel e_{w}^{j}(t)\parallel\leq\parallel e_{w}^{j}(t_{n})\parallel+C\Delta t\gamma^{2(j-1)}\parallel A_{22}^{-1}P_{1}^{*}(P_{2}^{*}\partial_{t}e_{w}^{0}(t)+\Pi_{1}A_{22}e_{w}^{0}(t))\parallel\nonumber\\
+&C\Delta t\frac{1}{1-\gamma^{2}}\left(\max_{0\leq i\leq j-1}\sup_{t\in(t_{n},t_{n+1})}\parallel e_{w}^{i}(t)\parallel+\max_{0\leq i\leq j-1}\sup_{t\in(t_{n},t_{n+1})}\parallel e_{u}^{i}(t)\parallel\right).
\end{align}
We set $\mu=C\Delta t\frac{1}{1-\gamma^{2}}\leq1$, then add \eqref{WR31} and \eqref{WR32}. With the condition $\parallel A_{11}^{-1}P_{2}^{*}A_{22}\parallel\leq C$ and $\parallel A_{22}^{-1}P_{1}^{*}A_{11}\parallel\leq C$, we can get 
\begin{align*}
&\sup_{t\in(t_{n},t_{n+1})}\parallel e_{u}^{j}(t)\parallel+\sup_{t\in(t_{n},t_{n+1})}\parallel e_{w}^{j}(t)\parallel\\
\leq&\mu\max_{0\leq i\leq j-1}\left(\sup_{t\in(t_{n},t_{n+1})}\parallel e_{u}^{i}(t)\parallel+\sup_{t\in(t_{n},t_{n+1})}\parallel e_{w}^{i}(t)\parallel\right)\\
&+C\sum_{i=0}^{1}\Delta t\mu^{i}\gamma^{2(j-1-i)}\sup_{t\in(t_{n},t_{n+1})}\parallel A_{22}^{-1}\partial_{t}e_{w}^{0}(t)+e_{w}^{0}(t)\parallel\\
\leq&C \sum_{i=0}^{j-1}\Delta t\mu^{i}\gamma^{2(j-1-i)}\sup_{t\in(t_{n},t_{n+1})}\parallel A_{22}^{-1}\partial_{t}e_{w}^{0}(t)+e_{w}^{0}(t)\parallel,
\end{align*}
for $n=1,2,\cdots,N-1$.
Then we get the desired result. $\Box$
%Finally we consider $\Delta t\rightarrow0$, thus $\mu\rightarrow0$, then we get the desired result. $\Box$

% We remark that, the convergence of the vanilla all-at-once algorithm based on the WR method for our partially explicit scheme is tested. Numerical results show that errors are decaying rapidly as the number of iterations increases. For example, takin g $\alpha=0.5$, $T=0.001$ and $\Delta t=10^{-5}$, i.e., the number of time steps is $100$, the method achieves $10^{-13}$ accuracy with only a few iteration steps (around $12$).

\subsection{Parareal all-at-once partially explicit temporal splitting algorithm}
\label{section:Parareal all-at-once partially explicit temporal splitting algorithm}
In this subsection, we describe in detail the use of the parareal algorithm to solve the partially explicit temporal splitting scheme and give its specific algorithm.

We further divide each time slice $[t_{n},t_{n+1}]$ into $M$ subintervals with $\delta t=\Delta t/P$, $\widehat{t}_p=p\delta t$, $p=0,1,\cdots,P$. According to the framework of the parareal algorithm, we should first compute a series of initial solutions based on the coarse solver $\mathcal{G}$. To do this, we solve the 
\eqref{eqnAd}-\eqref{eqnAe} sequentially with a coarse time step size to obtain the initial solutions. Then on each time slice $[t_{n},t_{n+1}]$,$n=0,1,\cdots,N-1$, we use the all-at-once method to solve \eqref{WR8}-\eqref{WR9}. This process iterates till convergence, which leads to our main algorithm: the parareal all-at-once partially explicit temporal splitting algorithm. 

% \begin{breakablealgorithm}
\begin{algorithm}
\caption{Parareal all-at-once partially explicit temporal splitting algorithm}
\begin{algorithmic}[1]
\Require initial date $u_{0}$, source term $f$, tolerance $\epsilon$, coarse matrices \eqref{eqnBh}, $\Psi_{1}$ and $\Psi_{2}$.
\Ensure $u_{H}^{n,k}$ for $n=1,2,\cdots,N$ .
\State Compute a series of initial solution $\{u_{n}^{0}\}_{n=1}^{N}$ and $\{w_{n}^{0}\}_{n=1}^{N}$
\begin{align*}
\binom{u_{n+1}^{0}}{w_{n+1}^{0}}=\mathcal{G}_{\Delta t}(\binom{u_{n}^{0}}{w_{n}^{0}},t_{n},t_{n+1}),\ n=0,1,\cdots,N-1.
\end{align*}
\For{$k=1,2,\cdots,$}
\State Parallel compute $\widehat{u}_{n+1}^{k}$ and $\widehat{w}_{n+1}^{k}$ on each time slice $[t_{n},t_{n+1}]$ 
\begin{align*}
\binom{\widehat{u}_{n+1}^{k}}{\widehat{w}_{n+1}^{k}}=\mathcal{F}_{\Delta t}(\binom{u_{n}^{k-1}}{w_{n}^{k-1}},t_{n},t_{n+1})
\end{align*}
\State Sequentially compute the corrected solution $\{u_{n}^{k+1}\}_{n=1}^{N}$ and $\{w_{n}^{k+1}\}_{n=1}^{N}$
\begin{align*}
\binom{u_{n+1}^{k+1}}{w_{n+1}^{k+1}}=\mathcal{G}_{\Delta t}(\binom{u_{n}^{k+1}}{w_{n}^{k+1}},t_{n},t_{n+1})+\binom{\widehat{u}_{n+1}^{k}}{\widehat{w}_{n+1}^{k}}-\mathcal{G}_{\Delta t}(\binom{u_{n}^{k}}{w_{n}^{k}},t_{n},t_{n+1})
\end{align*}
\State Determine whether the given condition is met.
\If $\max_{1\leq n\leq N}\parallel u_{n}^{k}-u_{n}^{k-1}\parallel<\epsilon$
\State return $u_{H}^{n,k}=u_{n}^{k}+w_{n}^{k},n=1,\cdots,N$
\State break
\EndIf
\EndFor
\end{algorithmic} \label{main_algorthm}
\end{algorithm}
% \end{breakablealgorithm}

We present in Algorithm \ref{main_algorthm} the main ingredient of the parareal all-at-once algorithm for \eqref{eqnAd}-\eqref{eqnAe}, which include the parallel computation in Step 3 and the sequential propagation in Step 4. 

\subsection{Convergence of the parareal algorithm}
\label{section:Error analysis}
This subsection is concerned with the convergence analysis for the Algorithm \ref{main_algorthm} in subsection \ref{section:Parareal all-at-once partially explicit temporal splitting algorithm}. To begin with, we introduce some lemmas that will be used in theoretical proof.
% \begin{lemma}\label{l1}
% Let $u$ the be solution of \eqref{eqnAa}, $u_{H}(t)$ be the solution of \eqref{eqnAb}, if $f_{t}\in L^{1}(0,T;L^{2}(\Omega))$ and $u_{tt}\in L^{1}(0,T;L^{2}(\Omega))$, then there holds\cite{reference8}
% \begin{align}\label{eqnCac}
% \begin{aligned}
% \parallel u(\cdot,T)-u_{H}(T)\parallel&\leq CH^{2}\Lambda^{-1}\kappa_{0}^{-\frac{1}{2}}\biggl(\max_{0\leq t\leq T}\parallel\kappa^{-\frac{1}{2}}(f-u_{t})\parallel\\
% &+\parallel\kappa^{-\frac{1}{2}}(f_{t}-u_{tt})\parallel_{L^{1}(0,T;L^{2}(\Omega))}\biggr)+\parallel u_{0}-u_{H,0}\parallel,
% \end{aligned}
% \end{align}
% where $\kappa_{0}=\max{\kappa^{-\frac{1}{2}}}$ and $\Lambda=\min_{1\leq j\leq N_{e}}\lambda_{J_{i}+1}^{(i)}$.
% \end{lemma}
% This lemma gives an error estimate in the spatial semi-discrete scheme.
\begin{lemma}\label{l2}
The coarse solve $\mathcal{G}$ is Lipschitz, it holds\cite{reference18}
\begin{align}\label{wr33}
\parallel\mathcal{G}_{\Delta t}\left(u,t_{n},t_{n+1}\right)-\mathcal{G}_{\Delta t}\left(v,t_{n},t_{n+1}\right)\parallel\leq\left(1+C_1\Delta t\right)\parallel u-v\parallel.
\end{align}
\end{lemma}

\begin{assumption}\label{l3}
Define $\mathcal{S}_{\Delta t}\left(u,t_{n},t_{n+1}\right)=\mathcal{F}_{\Delta t}\left(u,t_{n},t_{n+1}\right)-\mathcal{G}_{\Delta t}\left(u,t_{n},t_{n+1}\right)$. Then it has following property\cite{reference18}
\begin{align}\label{wr34}
\parallel\mathcal{S}_{\Delta t}\left(u,t_{n},t_{n+1}\right)\parallel\leq C_2(\Delta t)^{m+1}\parallel u\parallel.
\end{align}
\end{assumption}

% \textcolor{red}{
% \begin{lemma}\label{l4}
% The solution of the \eqref{eqnAb} is stable, in the sense that \cite{reference18}
% \begin{align*}
% \Vert u_H(t)\Vert\leq C\Vert u_{H,0}\Vert,
% \end{align*}
% where $C$ is independent of any parameters.
% \end{lemma}
% }
\textbf{Remark 1:}\quad Notice that in parareal algorithm, fine solver $\mathcal{F}$ is usually considered exact, so $m$ in Lemma \ref{l3} is the order of the coarse solver $\mathcal{G}$. Due to the fact that the coarse solver $\mathcal{G}$ is a first order scheme in our algorithm, thus we can safety replace $m$ by $1$ in the later analysis. 

\indent With the aid of above lemmas, with a similar proof as in \cite{reference42}, we have the following error results in full discrete scheme.
\begin{theorem} \label{thm:parareal}
Let $u_{H}(t_N)$ and $u_{H}^{N,k}$ be the solutions of \eqref{eqnAd}-\eqref{eqnAe} computed by fine solver and by Algorithm \ref{main_algorthm}, respectively, then we have
\begin{align*}
\parallel u_{H}(t_N)&-u_{H}^{N,k}\parallel\leq (C_2\Delta t T)^{k}\parallel u_{H,0}\parallel
\end{align*}
\end{theorem}
\textbf{Proof:} Following the framework of the Parareal algorithm, assume that the fine solver $\mathcal{F}$ is exact, i.e., $\forall n=0,1,\cdots,N-1$,
\begin{align*}
u_{H}(t_{n+1})=\mathcal{F}_{\Delta t}\left(u_{H}(t_{n}),t_{n},t_{n+1}\right).
\end{align*}
We have following equation
\begin{align*}
&u_{H}(T)\\
=&\mathcal{G}_{\Delta t}\left(u_{H}(t_{N-1}),t_{N-1},T\right)+\mathcal{F}_{\Delta t}\left(u_{H}(t_{N-1}),t_{N-1},T\right)-\mathcal{G}_{\Delta t}\left(u_{H}(t_{N-1}),t_{N-1},T\right)\\
=&\mathcal{G}_{\Delta t}\left(u_{H}(t_{N-1}),t_{N-1},T\right)+\mathcal{S}_{\Delta t}\left(u_{H}(t_{N-1}),t_{N-1},T\right).
\end{align*}
On the other hand, by \eqref{eqnAg} we have
\begin{align*}
u_{H}^{N,k}=&\mathcal{G}_{\Delta t}\left(u_{H}^{N-1,k},t_{N-1},T\right)+\mathcal{F}_{\Delta t}\left(u_{H}^{N-1,k-1},t_{N-1},T\right)-\mathcal{G}_{\Delta t}\left(u_{H}^{N-1,k-1},t_{N-1},T\right)\\
=&\mathcal{G}_{\Delta t}\left(u_{H}^{N-1,k},t_{N-1},T\right)+\mathcal{S}_{\Delta t}\left(u_{H}^{N-1,k-1},t_{N-1},T\right).
\end{align*}
Then by triangle inequality,
\begin{align}\label{WR35}
\begin{aligned}
\parallel u_{H}(T)-u_{H}^{N,k}\parallel=&\parallel(\mathcal{G}_{\Delta t}\left(u_{H}(t_{N-1}),t_{N-1},T\right)-\mathcal{G}_{\Delta t}\left(u_{H}^{N-1,k},t_{N-1},T)\right)\\
+&\left(\mathcal{S}_{\Delta t}(u_{H}(t_{N-1}),t_{N-1},T)-\mathcal{S}_{\Delta t}(u_{H}^{N-1,k-1},t_{N-1},T)\right).
\end{aligned}
\end{align}
With the help of the Lemma \ref{l2} and Assumption \ref{l3}, we get
\begin{align*} 
 \theta_{N}^k \leq (1+C_1 \Delta t)\theta_{N-1}^k+ C_2  \Delta t^2 \theta_{N-1}^{k-1},
\end{align*}
where $\theta_{N}^k = \parallel u_{H}(T)-u_{H}^{N,k}\parallel$, and $ \theta_{N}^0 \leq (1+C_1 \Delta t)\theta_{N-1}^0+ C_3  \Delta t^2$.
Multiply the equation $e_{n}^k = (1+C_1 \Delta t)e_{n-1}^k+ C_2  \Delta t^2 e_{n-1}^{k-1}$ by $\alpha^n$ and let $\rho^k(\alpha) = \sum_{n\geq 1} \alpha^n e_{n}^k$, we have
\begin{align*}
\begin{aligned}
    \rho^k(\alpha) & = (1+C_1 \Delta t) \alpha \rho^k(\alpha)+ C_2  \Delta t^2  \alpha \rho^{k-1}(\alpha),\\
    \rho^0(\alpha) &= (1+C_1 \Delta t) \alpha \rho^0(\alpha) + C_3  \Delta t^2 \frac{\alpha}{1-\alpha} 
\end{aligned}
\end{align*} 
By induction, we have 
\begin{align*}
     \rho^k(\alpha) =  C_3  \Delta t^2 ( C_2  \Delta t^2)^k 
    \frac{\alpha^{k+1}}{(1-\alpha)(1-(1+C_1\Delta t)\alpha)^{k+1}}.
\end{align*}
Then we have 
\begin{align*}
\begin{aligned}
  \rho^k(\alpha) &\leq C_3  \Delta t^2 ( C_2  \Delta t^2)^k 
    \frac{\alpha^{k+1}}{(1-(1+C_1\Delta t)\alpha)^{k+2}}\\
    & = C_3  \Delta t^2 ( C_2  \Delta t^2)^k 
    \alpha^{k+1} \sum_{n\geq 0}  \binom{n+k+1}{n} \alpha^n (1+C_1\Delta t)^n\\
    & \leq  \frac{C_3}{C_2}( C_2  \Delta t ^2)^{k+1} 
\sum_{m}  \binom{m}{m-k-1} \alpha^m (1+C_1\Delta t)^{m-k-1}
    \end{aligned}
\end{align*}
then 
\begin{align*}
\begin{aligned}e_{n}^k \leq  \frac{C_3}{C_2(k+1)!} (C_2\Delta t^2)^{k+1} e^{C_1(t_n - t_{k+1})} n^{k+1} \leq C (C_2\Delta t  t_n)^{k+1}     \end{aligned}
\end{align*}
That completes the proof.

\subsection{Parallel speedup analysis}
\label{section: Parallel speedup analysis}
Similar to the analysis for the parallel speedup and efficiency of the classical parareal algorithm presented in \cite{reference37}, we ignore the communication overhead, and adopt the following notation:
\begin{itemize}
\item{$T$ is the length of the time interval.}
\item{$N$ denote the number of coarse time interval as well as the number of processors.}
\item{$\Delta t$ is the time increment for the coarse propagator as well as the length of each processor.}
\item{$\delta t$ is the time increment for the fine propagation.}
\item{$K_0$ is the iteration number of the WR on each time sub-interval.}
\item{$N_1$ is the number of processors used in All-at-once method.}
\item{$K_1$ is the iteration number of the parareal all-at-once algorithm.}
\item{$S$ is the parallel speedup.}
\item{$E$ is the parallel efficiency with $E=S/N$.}
\end{itemize}

We assume that the time complexity for solving a $d_1-$dimensional linear system is $\tau_1$ and solving a $d_2-$dimensional linear system is $\tau_2$. 

The cost of computing the initial values with coarse propagation is $N\left(\tau_1+\tau_2\right)$. Since $P=\Delta t/\delta t$, then for the three steps in \eqref{WR7}, it has been proved in \cite{reference35} that step (a) and step (c) can be finished with the total cost $\mathcal{O}(d_1 P\log{(P)}/N_1)$, therefore the total cost of each iteration of the parareal all-at-once algorithm is $K_0\left(\mathcal{O}(d_1 P\log{(P)}/N_1)+P\tau_2\right)$. Then the total cost for the parareal all-at-once method is 
\begin{align*}
N\left(\tau_1+\tau_2\right)+K_1\left(K_0\left(Cd_1 P\log{(P)}/N_1+P\tau_2\right)+N\left(\tau_1+\tau_2\right)\right),
\end{align*}
where $C$ is a constant. The cost of applying $\mathcal{F}$ serially is $NP\left(\tau_1+\tau_2\right)$, hence the speedup for our proposed algorithm is 
\begin{align}\label{WR39}
\begin{aligned}
S=&\frac{NP\left(\tau_1+\tau_2\right)}{N\left(\tau_1+\tau_2\right)+K_1\left(K_0\left(Cd_1 P\log{(P)}/N_1+P\tau_2\right)+N\left(\tau_1+\tau_2\right)\right)}\\
=&\frac{NP\left(\tau_1+\tau_2\right)}{(1+K_1)N\left(\tau_1+\tau_2\right)+K_1K_0\left(Cd_1P\log{(P)}/N_1+P\tau_2\right)}\\
=&\frac{N}{(1+K_1)\frac{N}{P}+\frac{K_1K_0}{P\left(\tau_1+\tau_2\right)}\left(Cd_1 P\log{(P)}/N_1+P\tau_2\right)}.
\end{aligned}
\end{align}
Then the parallel efficiency is given by 
\begin{align}\label{WR40}
E=\frac{S}{N}=\frac{1}{(1+K_1)\frac{N}{P}+\frac{K_1K_0}{P\left(\tau_1+\tau_2\right)}\left(Cd_1 P\log{(P)}/N_1+P\tau_2\right)}.
\end{align}
If we consider the case $N=P$, i.e., $\frac{T}{\Delta t}=\frac{\Delta t}{\delta t}$, then the parallel efficiency is bounded by
\begin{align*}
E\leq\frac{1}{1+K_1}.
\end{align*}

\section{Numerical experiments}
\label{section:Numerical experiments}
In this section, we perform some numerical experiments to verify the feasibility and effectiveness of the Algorithm \ref{main_algorthm} and the fine propagator. In all numerical examples, we consider the computational domain $\Omega=[0,1]\times[0,1]$, the coarse scale and fine scale spatial mesh size are $H=\frac{1}{10}$ and $h=\frac{1}{100}$, respectively. The relative error is defined as follows:
\begin{align*}
\frac{\Vert u_h^N-u_H^{N}\Vert_{L^{2}(\Omega)}}{\Vert u_h^N\Vert_{L^2(\Omega)}}
\end{align*}
where $u_h$ denotes the reference solution obtained by the finite element method in space and the WR method in time. Moreover, we use $R_p$ and $R_w$ represent the theoretical rate of the Algorithm \ref{main_algorthm} and WR method in the semilog plots, respectively.
\subsection{Numerical experiment 1}
\label{section: Numerical experiment 1}
\begin{figure}[H]
\centering\begin{minipage}{0.4\linewidth}
\includegraphics[width=2.05in,height=1.7in]{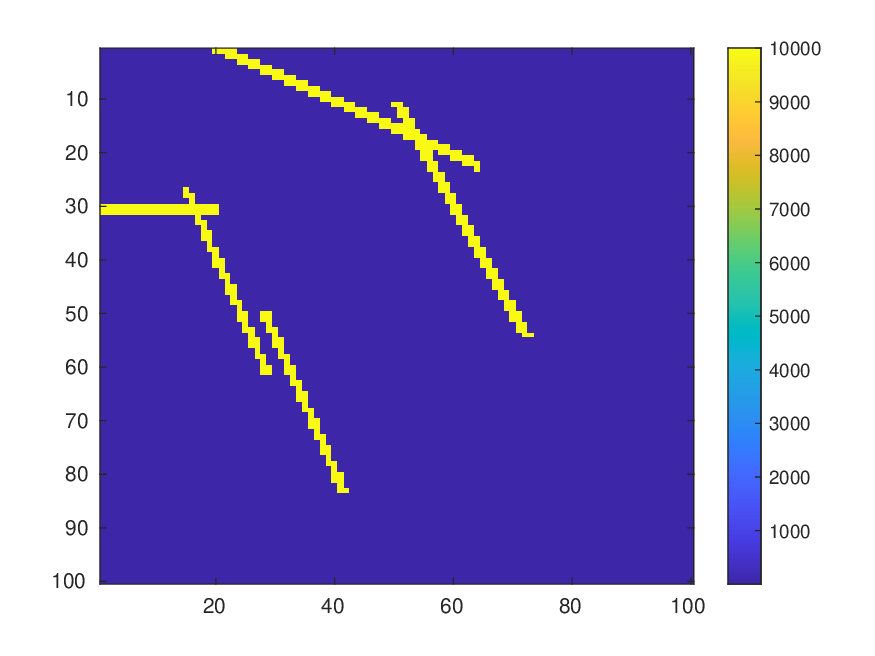}
\end{minipage}
\centering\begin{minipage}{0.4\linewidth}
\includegraphics[width=2.05in,height=1.7in]{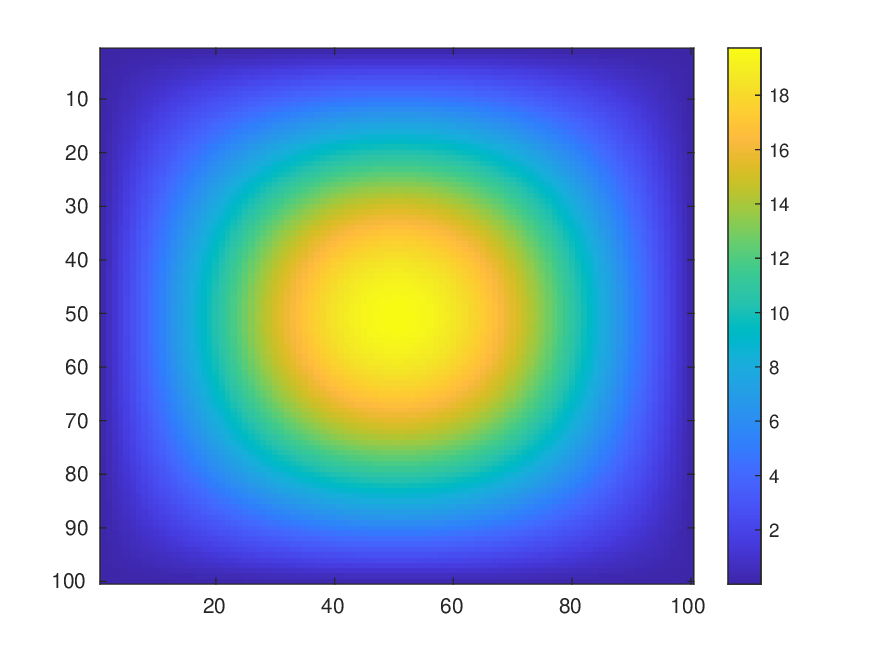}
\end{minipage}
\caption{Example 1. \textbf{Left:} $\kappa$; \textbf{Right:} $f$}\label{f2}
\end{figure}

The medium parameter $\kappa$ and the source term $f$ are shown in Figure \ref{f2}. We see that the permeability field is heterogeneous with high contrast. The contrast is $10^4$. In this example, we consider the total simulation time $T=0.001$, and set the tolerance $\epsilon=10^{-13}$.

\begin{figure}[H]
\centering\begin{minipage}{0.4\linewidth}
\includegraphics[width=2.05in,height=1.7in]{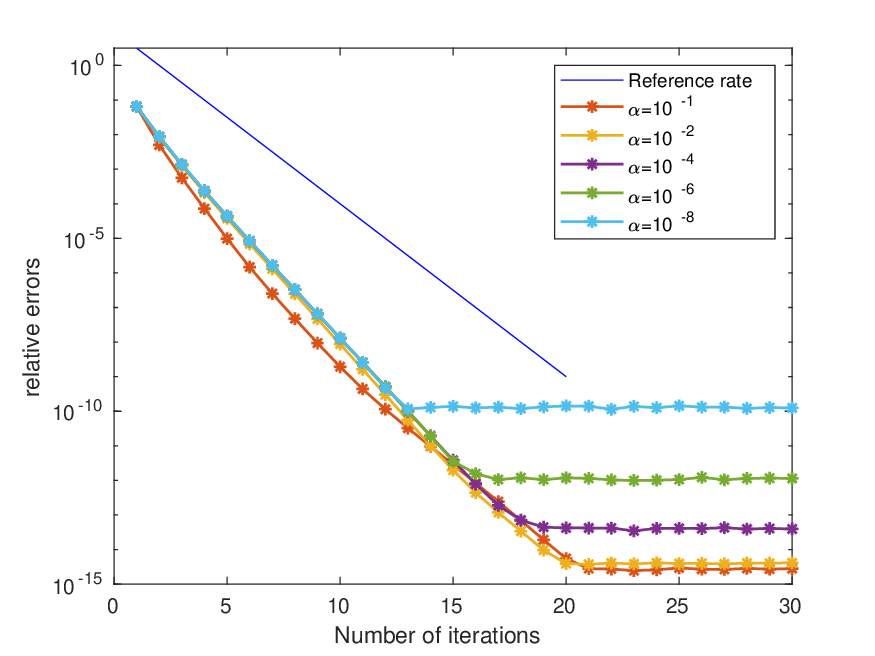}
\end{minipage}
\centering\begin{minipage}{0.4\linewidth}
\includegraphics[width=2.05in,height=1.7in]{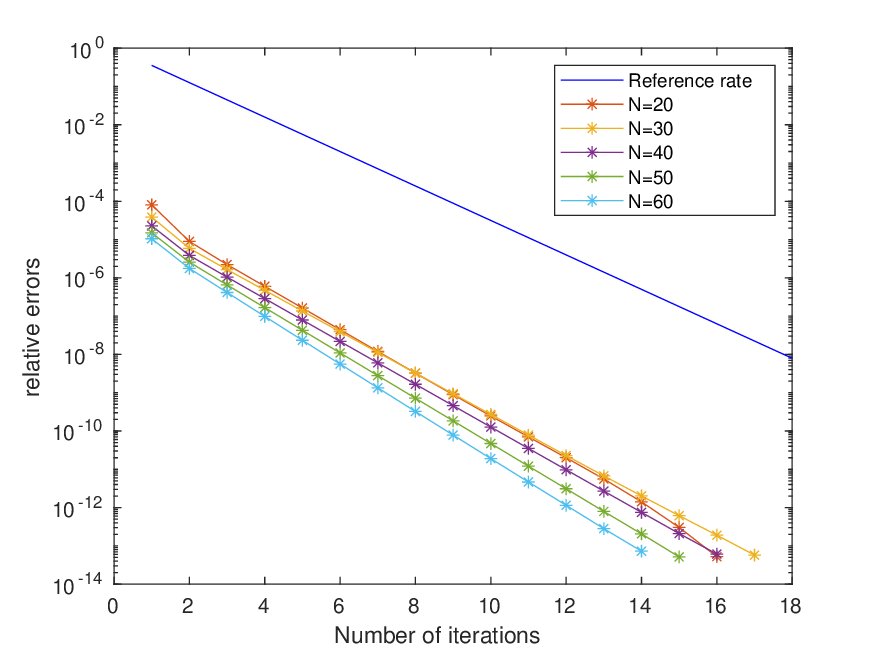}
\end{minipage}
\caption{Example 1. \textbf{Left:} Convergence rate of the WR method with $R_w=-0.56$; \textbf{Right:} Convergence rate of the Algorithm \ref{main_algorthm} with $R_p=-0.5$.}\label{f4}
\end{figure}

We first investigate the influence of the parameter $\alpha$ on the convergence of the WR method. Remark that we adopt the WR method at each subinterval $(t_n,t_{n+1})$, then we only verify the convergence result at interval $(t_0,t_1)$. According to Theorem 4.1, the theoretical rate for the log of error against the number of iterations is $log(\gamma^2)$ where $\gamma = \sup_{v_1 \in V_{H,1}, v_2 \in V_{H,2}}\frac{|(v_{1},v_{2})|}{\parallel v_{1}\parallel\parallel v_{2}\parallel}$. We numerically compute the value of $\gamma$. We choose $N=20$, $P=30$ with $\delta t=1/6\times10^{-5}$ and $\alpha=10^{-1}, 10^{-2}, 10^{-4}, 10^{-6}, 10^{-8}$ respectively. In the left of Figure \ref{f4}, we present the convergence behavior of the WR method for different values of $\alpha$, along with its theoretical convergence rate. It is evident that a larger value of $\alpha$ leads to slightly faster convergence, albeit at the cost of a relatively lower accuracy- a finding consistent with the results reported in \cite{reference35}.

\begin{figure}[H]
\centering\begin{minipage}{0.4\linewidth}
\includegraphics[width=2.05in,height=1.7in]{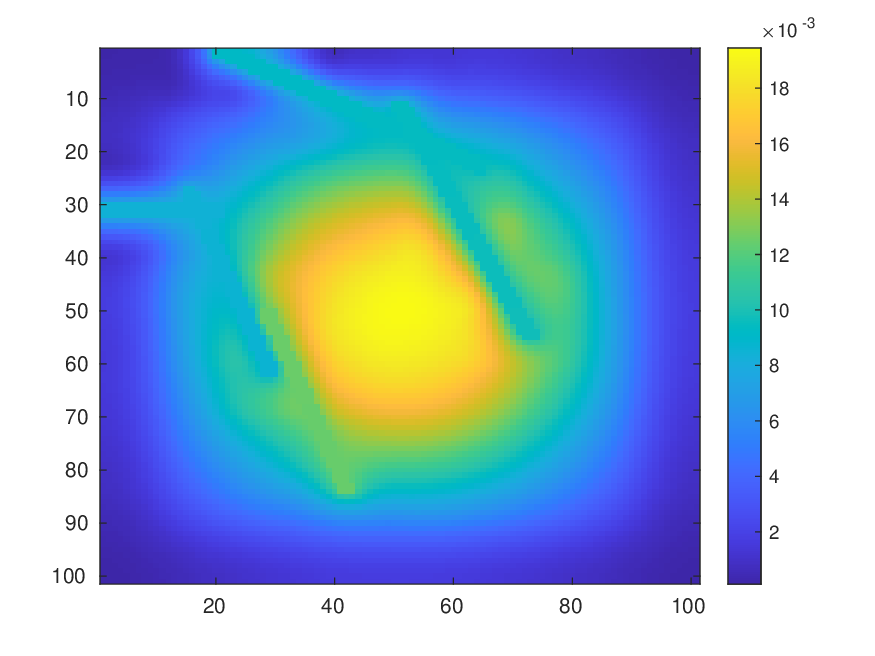}
\end{minipage}
\centering\begin{minipage}{0.4\linewidth}
\includegraphics[width=2.05in,height=1.7in]{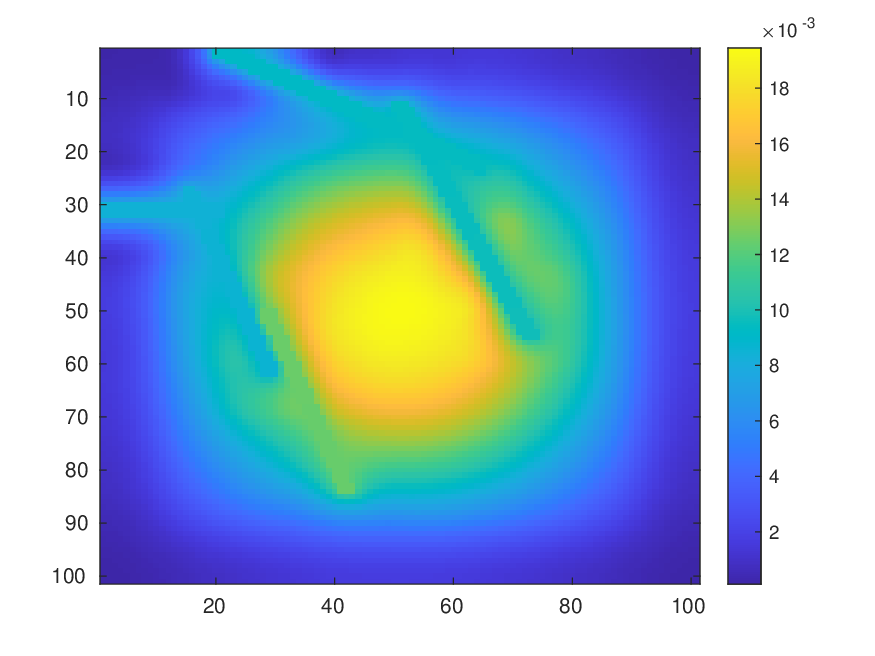}
\end{minipage}
\caption{Example 1. \textbf{Left:} reference solution; \textbf{Right:} Algorithm \ref{main_algorthm} solution with $R_p=-0.5$.}\label{f3}

\end{figure}

We then investigate the feasibility and effectiveness of the Algorithm \ref{main_algorthm}. We choose $\alpha=0.1$, $N=20,30,40,50,60$ and keep $NP=600$ and $\delta t=1/6\times10^{-5}$. The reference solution and the solution obtained by Algorithm \ref{main_algorthm} at the final time with $N=50$ are presented in Figure \ref{f3}.  In the right of Figure \ref{f4}, we show the convergence behavior of the Algorithm \ref{main_algorthm} for different values of $N$, along with its theoretical convergence rate. According to Theorem 4.2, the theoretical rate for the log of error against the number of iterations is $log(C_2 \Delta t T)$. We numerically calculate this value and it is around or greater than $-0.5 \approx log(0.316)$ in examples \cite{reference43}, so we uniformly take $R_p=-0.5$ as the theoretical rate in the semilog plots. It can be observed that our proposed algorithm shows fast convergence.

\begin{figure}[H]
\centering\begin{minipage}{0.4\linewidth}
\includegraphics[width=2.05in,height=1.7in]{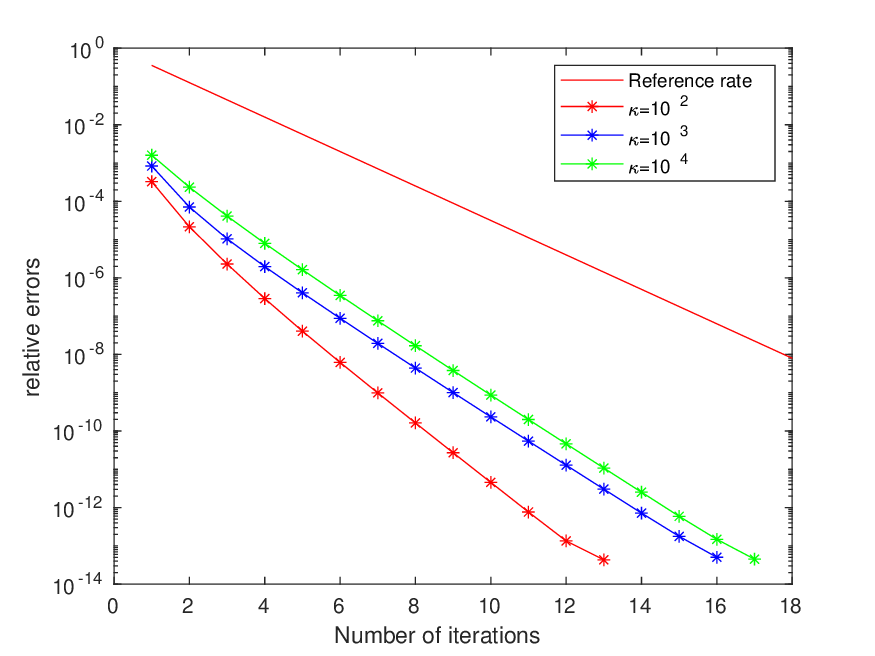}
\end{minipage}
\centering\begin{minipage}{0.4\linewidth}
\includegraphics[width=2.05in,height=1.7in]{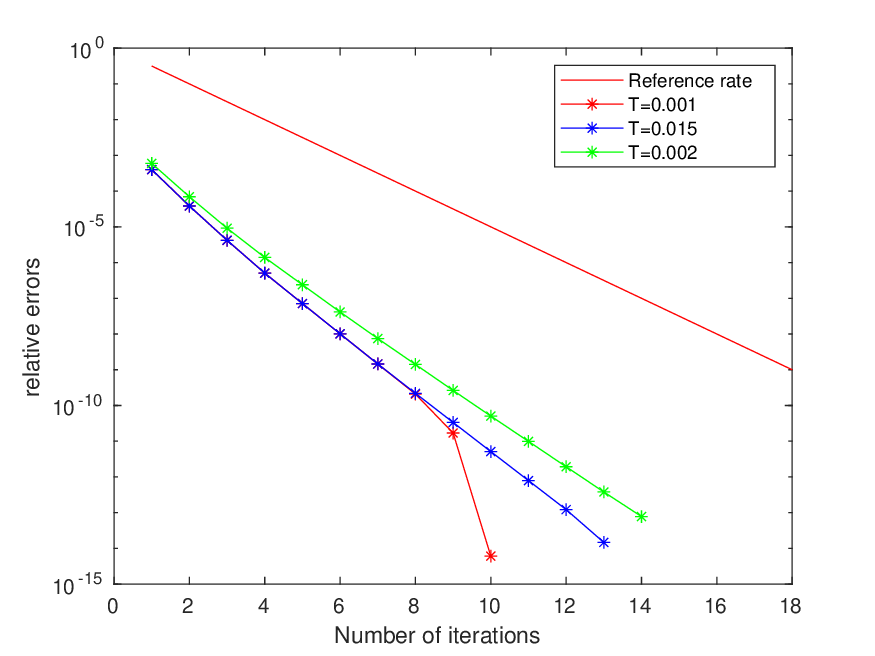}
\end{minipage}
\caption{Example 1. \textbf{Left:} Convergence rate of the Algorithm \ref{main_algorthm} for different contrast ratio; \textbf{Right:} Convergence rate of the Algorithm \ref{main_algorthm} for different $T$. $R_p=-0.5$.}\label{ff1}
\end{figure}

Finally, we investigate the influence of the high contrast ratio and the parameter $T$ on the convergence rate. In the left of Figure \ref{ff1}, we choose $T=0.001$, $N=50$, $NP=600$ with $\delta t=1/6\times10^{-5}$, and the contrast ratios are considered as $10^2$, $10^3$ and $10^4$ respectively. Clearly, a larger contrast in $\kappa$ leads to a larger error, which is accompanied by a slightly slower convergence rate. In the right of Figure \ref{ff1}, we fixed $\Delta t=1\times10^{-4}$ and $\delta t=5\times10^{-6}$ to investigate the influence of $T$ on the convergence rate. We consider $T=0.001,0.0015,0.002$, $N=10,15,20$, $NP=200,300,400$, respectively. We take $L^2$ errors in space and $L^{\infty}$ error in time. We see that as $T$ increases, the convergence become a little bit slower, which agrees with the theoretical result. We remark that, the reason for the sharp drop in the red error line at the final step in right of Figure \ref{ff1} is that the numerical solution has converged to the reference solution at the $10$-th iteration (since the number of coarse time step size is $10$ in this case).

% In the right of Figure \ref{f4}, we give the max differences computed by \eqref{eqnAh} for different $N$. It is clear that as the number of iterations increase, the max differences become smaller and eventually stabilize. It can be observed that as the number of processors (i.e., $N$) increases, the number of iterations needed to reach the desired accuracy decreases and eventually stabilizes.


\subsection{Numerical experiment 2}
\begin{figure}[H]
\centering\begin{minipage}{0.4\linewidth}
\includegraphics[width=2.05in,height=1.7in]{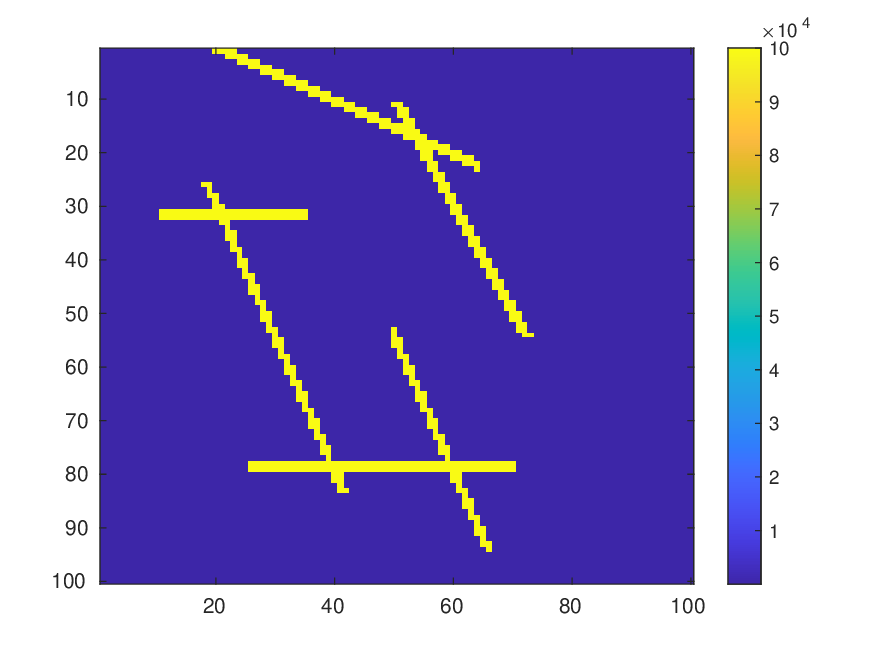}
\end{minipage}
\centering\begin{minipage}{0.4\linewidth}
\includegraphics[width=2.05in,height=1.7in]{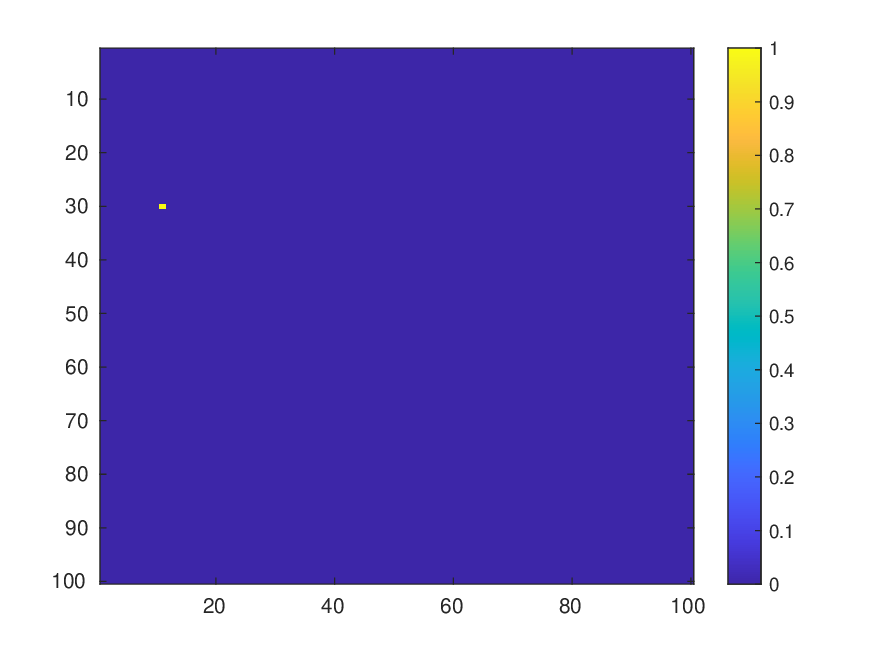}
\end{minipage}
\caption{Example 2. \textbf{Left:} $\kappa$; \textbf{Right:} $f$}\label{f6}
\end{figure}

In this numerical example, the multiscale high contrast coefficient $\kappa$ and the source term $f$ are presented in Figure \ref{f6}. The source term now is considered to be a point source term. We consider $T=0.01$ and set the tolerance $\epsilon=10^{-12}$.

\begin{figure}[H]
\centering\begin{minipage}{0.4\linewidth}
\includegraphics[width=2.05in,height=1.7in]{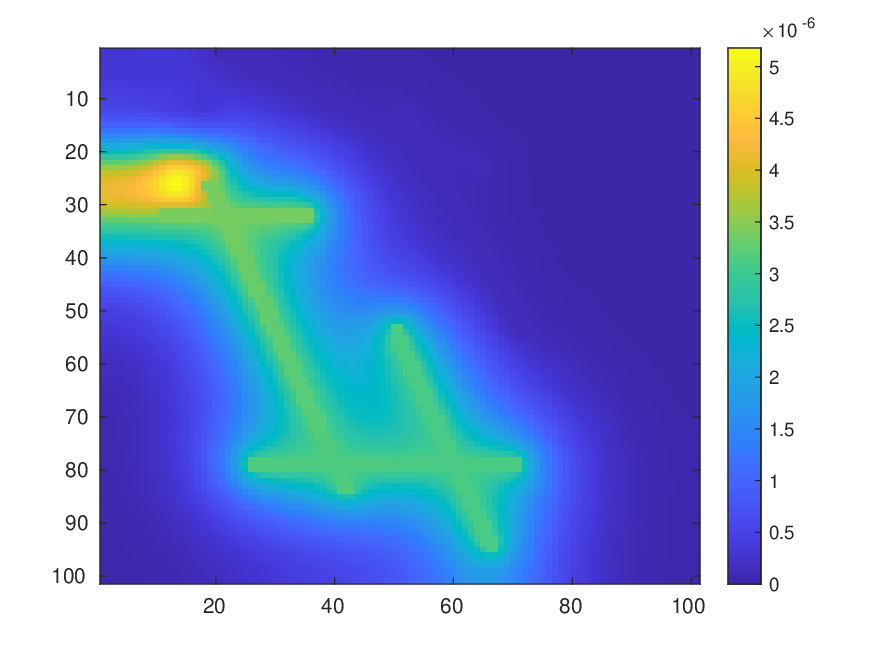}
\end{minipage}
\centering\begin{minipage}{0.4\linewidth}
\includegraphics[width=2.05in,height=1.7in]{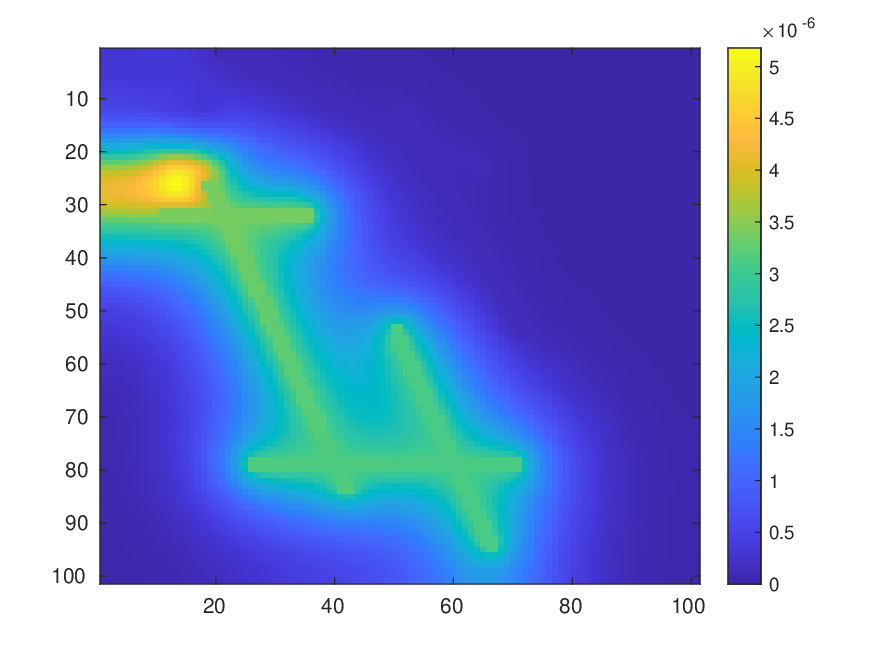}
\end{minipage}
\caption{Example 2. \textbf{Left:} reference solution; \textbf{Right:} Algorithm \ref{main_algorthm} solution.}\label{f7}
\end{figure}

\begin{figure}[H]
\centering\begin{minipage}{0.4\linewidth}
\includegraphics[width=2.05in,height=1.7in]{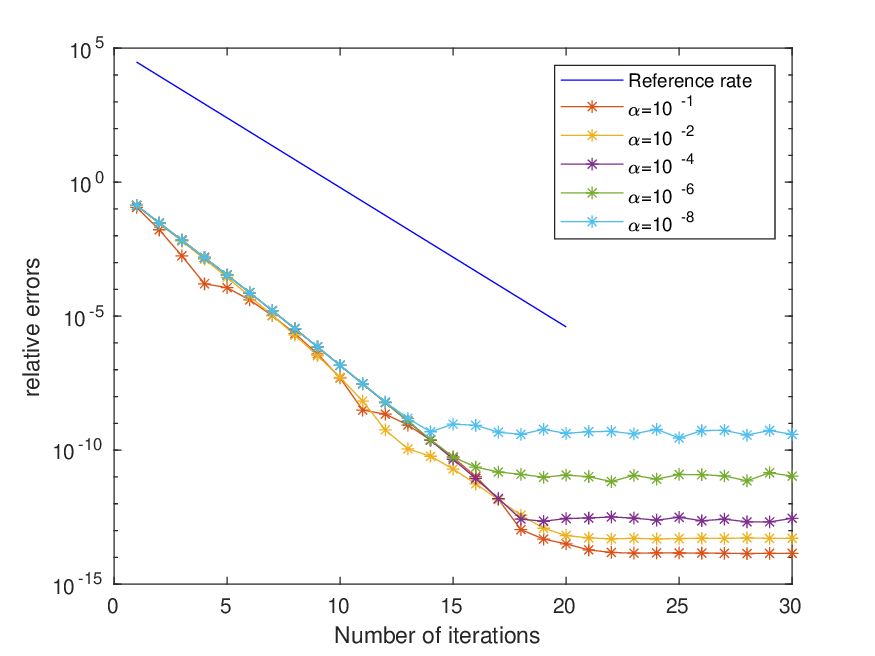}
\end{minipage}
\centering\begin{minipage}{0.4\linewidth}
\includegraphics[width=2.05in,height=1.7in]{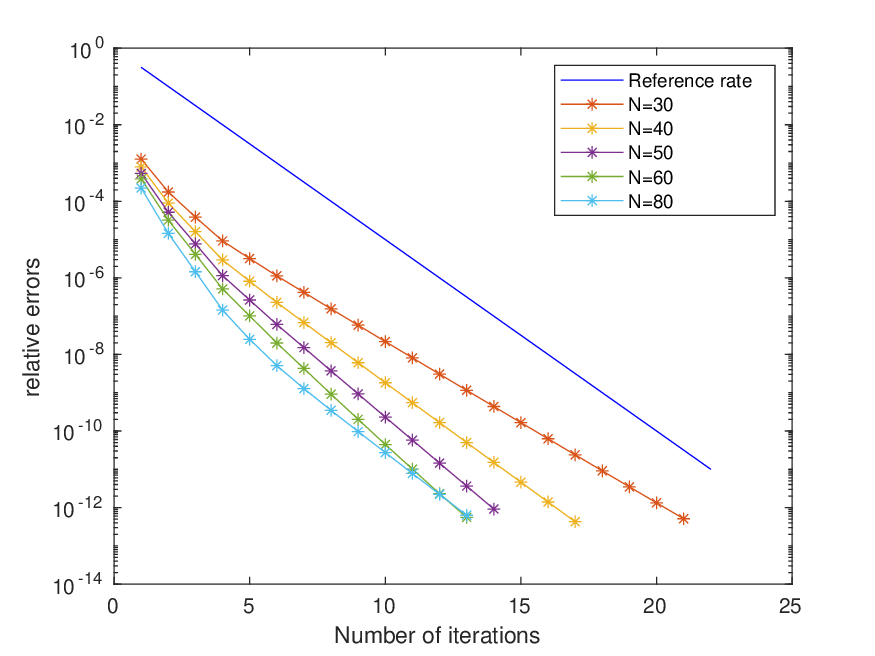}
\end{minipage}
\caption{Example 2. \textbf{Left:} Convergence rate of the WR technique with $R_w=-0.58$; \textbf{Right:} Convergence rate of the Algorithm \ref{main_algorthm} with $R_p=-0.5$.}\label{f8}
\end{figure}

\begin{figure}[H]
\centering\begin{minipage}{0.4\linewidth}
\includegraphics[width=2.05in,height=1.7in]{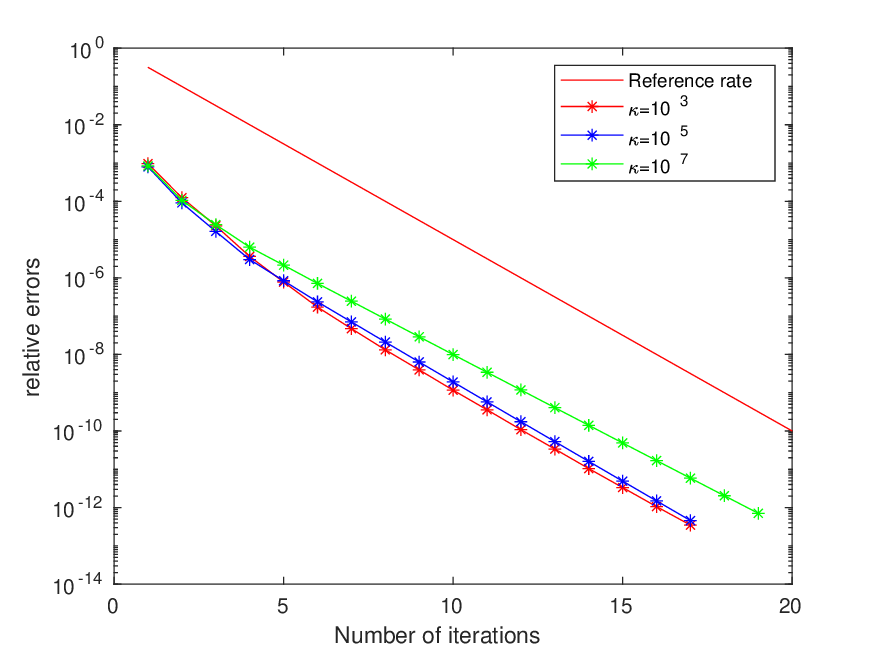}
\end{minipage}
\centering\begin{minipage}{0.4\linewidth}
\includegraphics[width=2.05in,height=1.7in]{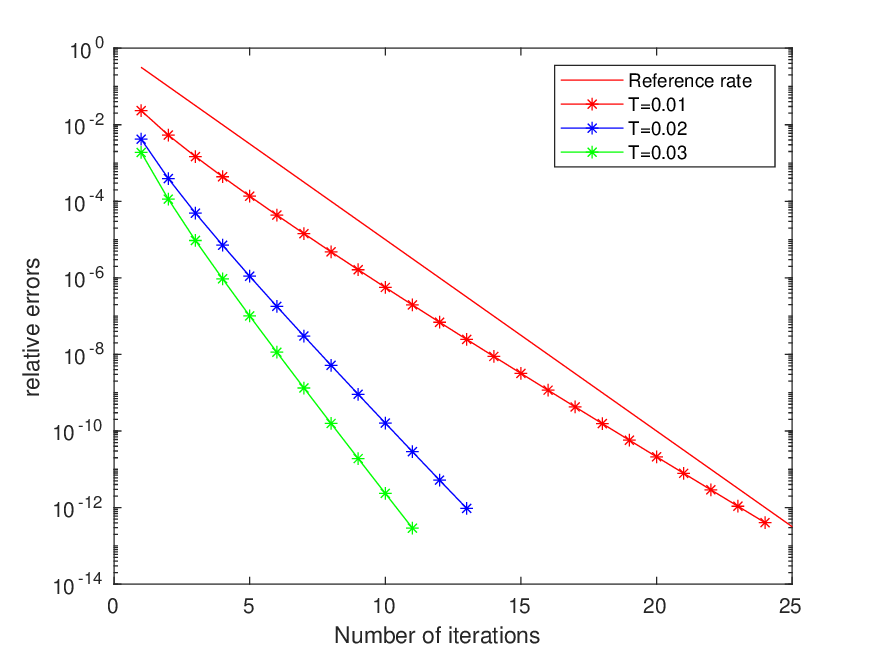}
\end{minipage}
\caption{Example 2. \textbf{Left:} Convergence rate of the Algorithm \ref{main_algorthm} for different contrast ratio; \textbf{Right:} Convergence rate of the Algorithm \ref{main_algorthm} for different $T$. $R_p=-0.5$.}\label{ff2}
\end{figure}

We first investigate the influence of the parameter $\alpha$ on the convergence of the WR method. We still verify the convergence behavior at $(t_0,t_1)$. Here we choose $N=30$, $P=40$ with $\delta t=25/3\times10^{-6}$, and $\alpha=10^{-1}, 10^{-2}, 10^{-4}, 10^{-6}, 10^{-8}$ respectively. We give the convergence result of the WR method for different values of $\alpha$ and its theoretical convergence rate in the left of Figure \ref{f8}. It is obvious that the parameter $\alpha$ influences both the convergence rate and the accuracy of the WR method. Then we investigate the behavior of the Algorithm \ref{main_algorthm}. We choose $\alpha=0.1$, $N=30,40,50,60,80$ and keep $NP=1200$ and $\delta t=25/3\times10^{-6}$. Figure \ref{f7} shows the reference solution and Algorithm \ref{main_algorthm} solution at $t=T$ with $N=50$. Similarly as in Example 1, we give the convergence rate along with its theoretical convergence rate in the right of Figure \ref{f8}. We can see that the error decays more rapidly as $N$ increases. In the left of Figure \ref{ff2}, we present the convergence rate for different contrast ratios with $T=0.01$, $N=40$, $NP=1200$ and $\delta t=25/3\times10^{-6}$. Similarly to experiment 1, the contrast ratio slightly affects the accuracy of the error and the convergence rate. In the right of Figure \ref{ff2}, we fix $\Delta t=1/3\times10^{-3}$ and $\delta t=1/6\times10^{-4}$ and consider $T=0.01,0.02,0.03$, $N=30,60,90$ and $NP=600,1200,1800$ respectively. We see that as $T$ increases, the log of convergence factor become smaller.

% \begin{figure}[H]
%     \centering
%     \includegraphics[width=0.6\linewidth]{figure/N2/WR_test2.eps}
%     \caption{Example 2. .}
%     \label{f9}
% \end{figure}
\subsection{Numerical experiment 3}
\begin{figure}[H]
\centering\begin{minipage}{0.4\linewidth}
\includegraphics[width=2.05in,height=1.7in]{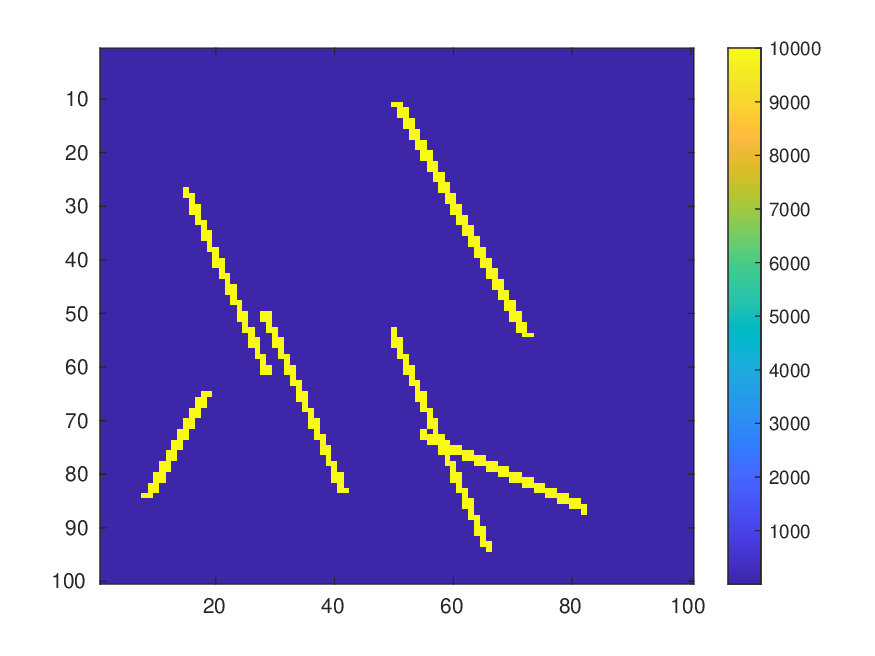}
\end{minipage}
\centering\begin{minipage}{0.4\linewidth}
\includegraphics[width=2.05in,height=1.7in]{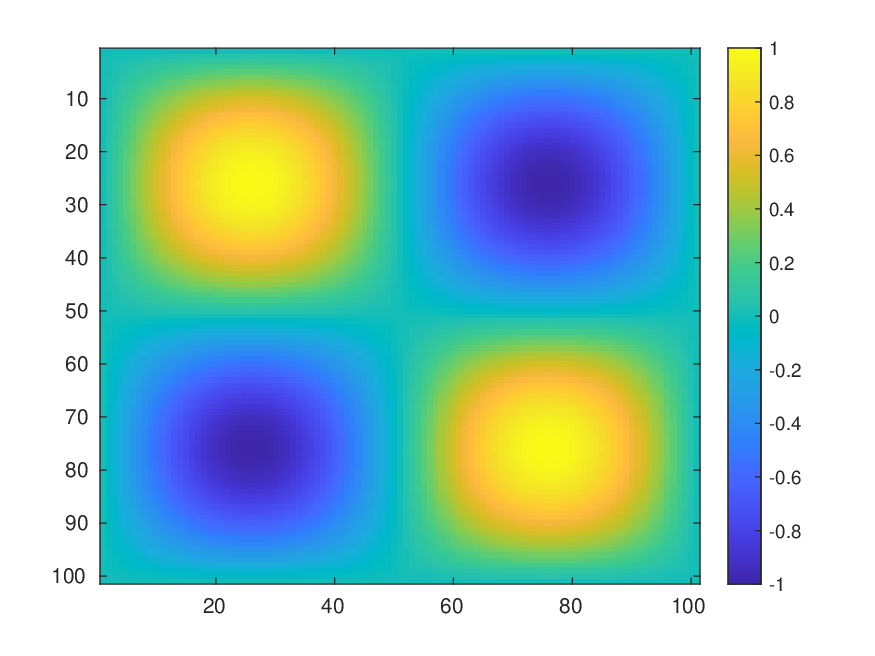}
\end{minipage}
\caption{Example 3. \textbf{Left:} $\kappa$; \textbf{Right:} $u_0(x,y)$.}\label{f10}
\end{figure}

In our final numerical experiment, we consider a more complicated multiscale parameter $\kappa$ and take $f=0$. The initial condition is taken to be $u_0(x,y)=\sin{(2\pi x)}\sin{(2\pi y)}$. We plot the multiscale coefficient $\kappa$ and the initial condition $u_0(x,y)$ in Figure \ref{f10}. The total simulation time $T=0.02$ and the tolerance is chosen to be $\epsilon=10^{-10}$.

\begin{figure}[H]
\centering\begin{minipage}{0.4\linewidth}
\includegraphics[width=2.05in,height=1.7in]{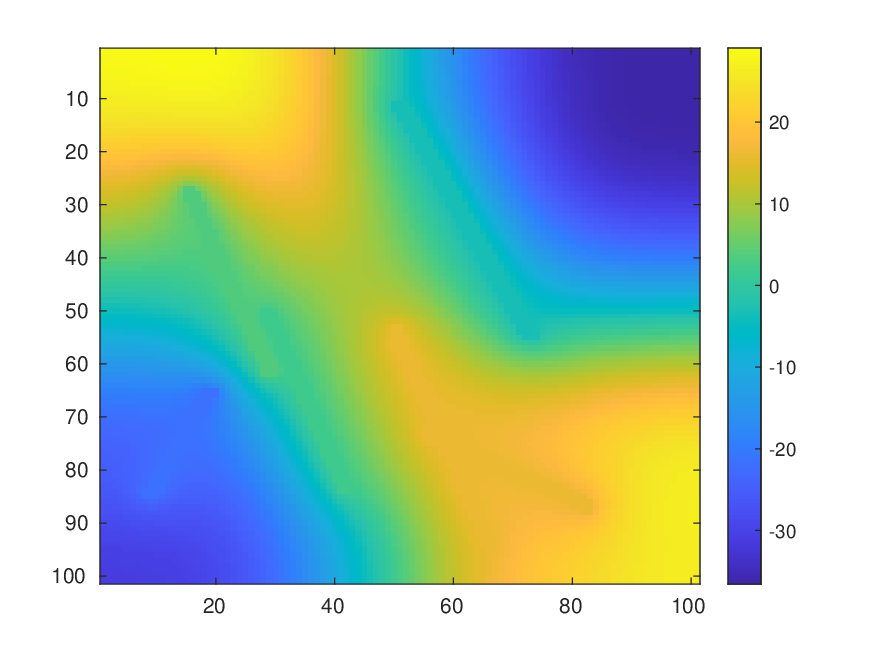}
\end{minipage}
\centering\begin{minipage}{0.4\linewidth}
\includegraphics[width=2.05in,height=1.7in]{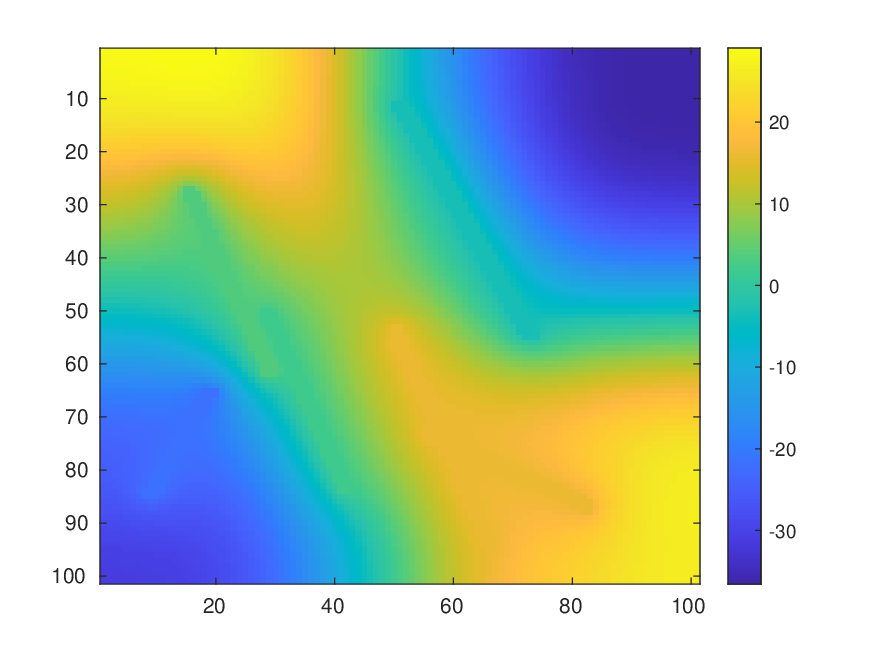}
\end{minipage}
\caption{Example 3. \textbf{Left:} reference solution; \textbf{Right:} Algorithm \ref{main_algorthm} solution.}\label{f11}
\end{figure}

\begin{figure}[H]
\centering\begin{minipage}{0.4\linewidth}
\includegraphics[width=2.05in,height=1.7in]{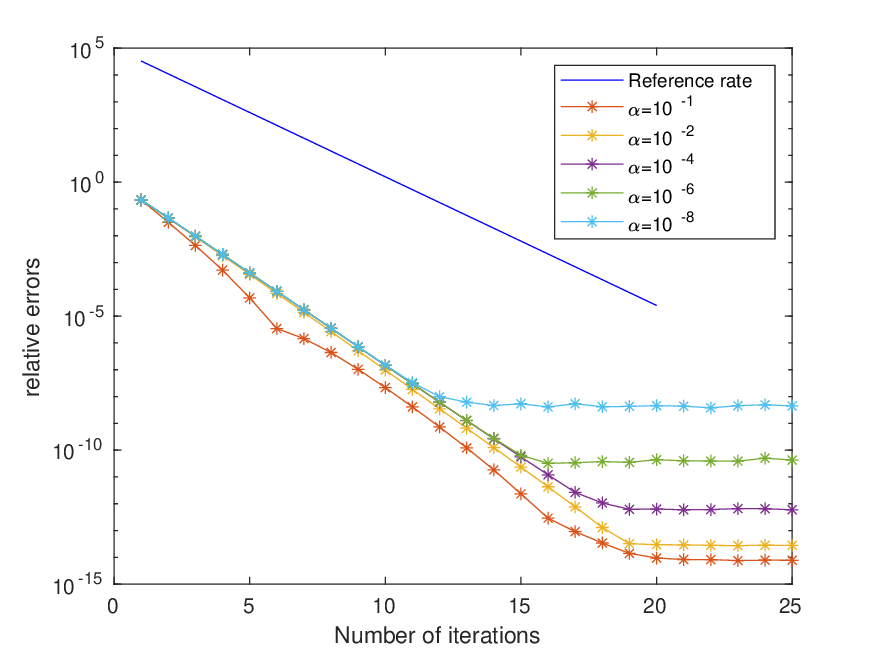}
\end{minipage}
\centering\begin{minipage}{0.4\linewidth}
\includegraphics[width=2.05in,height=1.7in]{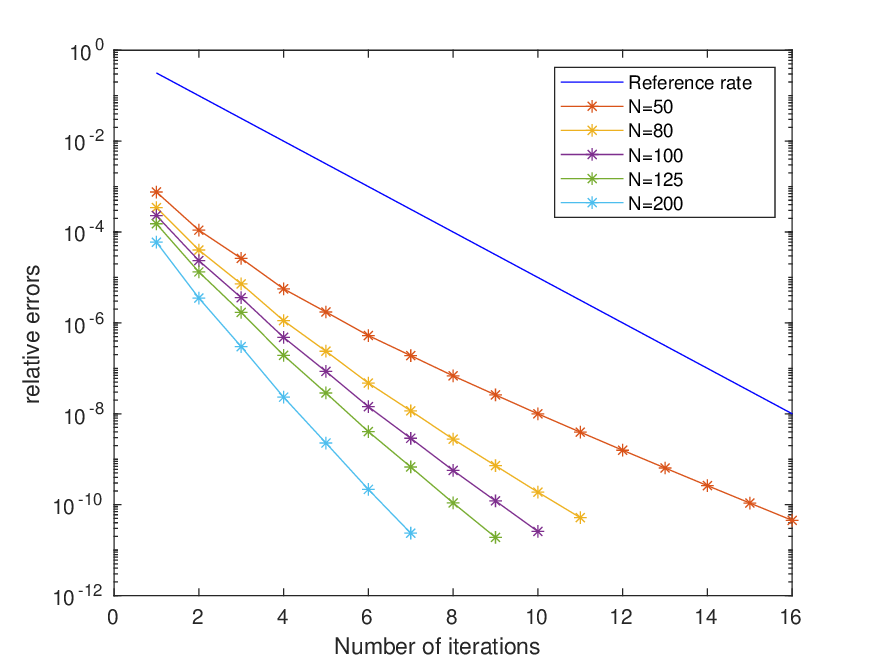}
\end{minipage}
\caption{Example 3. \textbf{Left:} Convergence rate of the WR technique with $R_w=-0.48$; \textbf{Right:} Convergence rate of the Algorithm \ref{main_algorthm} with $R_p=-0.5$.}\label{f12}
\end{figure}

\begin{figure}[H]
\centering\begin{minipage}{0.4\linewidth}
\includegraphics[width=2.05in,height=1.7in]{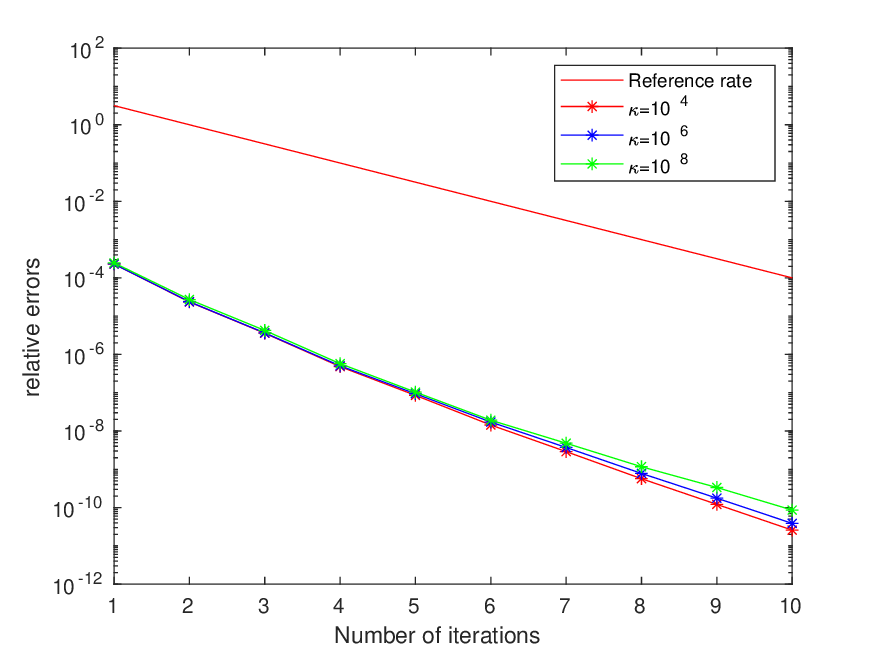}
\end{minipage}
\centering\begin{minipage}{0.4\linewidth}
\includegraphics[width=2.05in,height=1.7in]{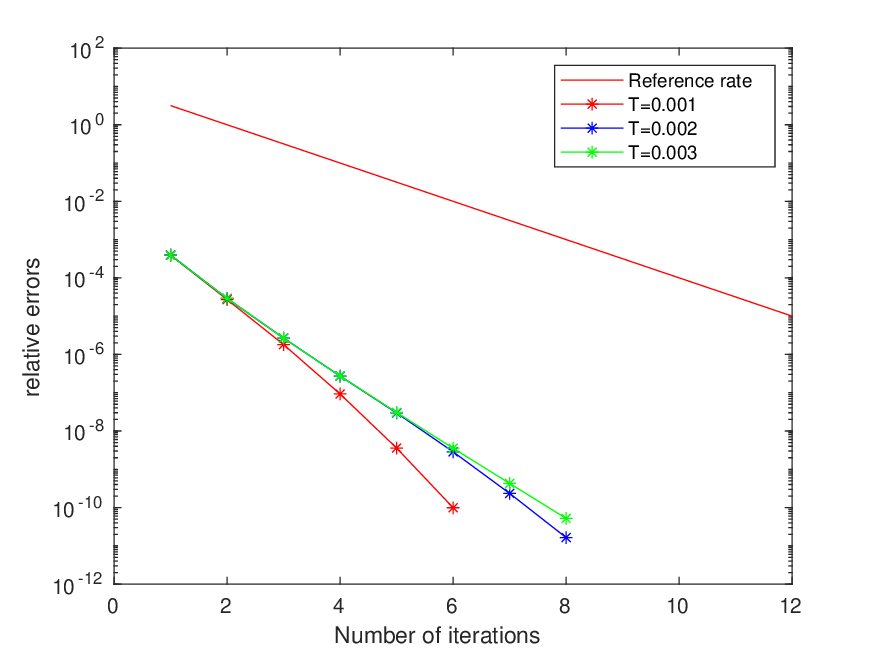}
\end{minipage}
\caption{Example 3. \textbf{Left:} Convergence rate of the Algorithm \ref{main_algorthm} for different contrast ratio; \textbf{Right:} Convergence rate of the Algorithm \ref{main_algorthm} for different $T$. $R_p=-0.5$.}\label{ff3}
\end{figure}

The left of Figure \ref{f12} provides the convergence rate of the WR method for different values of $\alpha$ and its theoretical convergence rate in $(t_0,t_1)$. Then we investigate the convergence of the Algorithm \ref{main_algorthm}. We choose $\alpha=0.1$, $N=50,80,100,125,200$ and keep $NP=2000$ and $\delta t=1\times10^{-5}$. The reference solution and the Algorithm \ref{main_algorthm} solution at $t=T$ with $N=200$ are shown in Figure \ref{f11}. The convergence rate along with its theoretical convergence rate is presented in the right of Figure \ref{f12}. In Figure \ref{ff3}, we present the convergence rate for different values of $\kappa$ with $T=0.02$, $N=100$ and $NP=2000$ with $\delta t=1\times10^{-5}$ and the convergence rate for $T=0.001,0.002,0.003$ with fixed $\Delta t=1\times10^{-4}$ and $\delta t=2\times10^{-6}$, respectively. Similar conclusion can be found in this example.

\textbf{Remark 2:} The theoretical rate in Theorem \ref{thm:parareal} was obtained based on the the assumption that the fine propagator is exact and coarse propagator is strongly stable, our numerical tests show that the convergence rate is a little bit slower for small $N$ (see the right of Figure \ref{f8} and \ref{f12}), but the convergence could be faster for large $N$, this might occur because the assumptions are not met, better results require more investigation in our future studies.

% \begin{figure}[H]
%     \centering
%     \includegraphics[width=0.6\linewidth]{figure/N3/WR_test3.eps}
%     \caption{Example 3. .}
%     \label{f13}
% \end{figure}

\section{Conclusion}
\label{section:Conclusion}
In this paper, we consider the diffusion equation with high contrast coefficient, and propose a parareal algorithm for the partially explicit temporal splitting scheme. In the parareal algorithm, we propose a one-step partially explicit temporal splitting scheme as the coarse solver, and utilize the all-at-once method as the fine solver to efficiently improve the computational efficiency. The convergences of the all-at-once method and the proposed parareal algorithm are given. An error estimate for the full discretization is given. Numerical experiments show that the proposed algorithm is computationally fast and accurate. The algorithm and the analysis in this paper are based on the case of a linear model. We will consider generalizing the algorithm to the nonlinear case in future work.

%\printbibliography
\bibliographystyle{plain} %plain
\bibliography{lookup}

\end{document}